\documentclass[12pt]{amsart}

\usepackage{latexsym, amssymb}

\newcommand{\be}{\begin{equation}}
\newcommand{\ee}{\end{equation}}
\newcommand{\beq}{\begin{eqnarray}}
\newcommand{\eeq}{\end{eqnarray}}

\newtheorem{thm}{Theorem}[section]
\newtheorem{prop}{Proposition}[section]
\newtheorem{cor}{Corollary}[section]
\newtheorem{lma}{Lemma}[section]

\newtheorem{rmk}{Remark}[section]

\def\Pi{\displaystyle{\mathbb{II}}}

\def\p{\partial}
\def\S{\Sigma}

\def\R{\mathbb{R}}
\def\H{\mathbb{H}}

\def\tr{{\rm tr}}
\def\p{\partial}
\def\lf{\left}
\def\ri{\right}
\def\e{\epsilon}
\def\ol{\overline}
\def\R{\Bbb R}

\def\la{\langle}
\def\ra{\rangle}
\def\hg{\hat g}
\def\hh{\hat h}

\def\hg{\hat{g}}
\def\hh{\hat{h}}
\def\l{\lambda}
\def\Ric{\text{\rm Ric}}
\def\div{\text{\rm div}}
\def\grad{\text{\rm grad}}

\def\hg{{\hat g}}

\def\gst{g_{ \mathbb{S}^{ m } } }

\def\hst{\nabla^2_{\gst}}

\begin{document}

\title{On the volume functional of compact manifolds
with boundary with constant scalar curvature}

\date{June 2008}

\author{Pengzi Miao and Luen-Fai Tam$^1$}

\renewcommand{\subjclassname}{%
  \textup{2000} Mathematics Subject Classification}

\thanks{$^1$Research partially supported by Earmarked Grant of Hong
Kong \#CUHK403005}

\address{The School of Mathematical Sciences, Monash University,
Victoria, 3800, Australia.} \email{Pengzi.Miao@sci.monash.edu.au}

\address{The Institute of Mathematical Sciences and Department of
 Mathematics, The Chinese University of Hong Kong,
Shatin, Hong Kong, China.} \email{lftam@math.cuhk.edu.hk}

\begin{abstract}
 We study the volume functional on the space of constant scalar curvature metrics 
 with a prescribed boundary metric. We derive a sufficient and necessary condition 
 for a metric to be
 a critical point, and show that the only domains in space forms, on which the standard
 metrics are critical points, are geodesic balls. In the zero scalar curvature case, assuming
 the boundary  can be isometrically embedded in the Euclidean space as a
 compact strictly convex hypersurface, we show that the volume of a critical point
 is always no less than the Euclidean volume bounded by the isometric embedding of
 the boundary, and the two volumes are equal if and only if the critical point is isometric to
 a standard Euclidean ball. We also derive a second variation formula and apply it to
 show that, on Euclidean balls and ``small'' hyperbolic and spherical balls in
 dimensions $3\le n\le 5$, the standard space form metrics are indeed saddle
 points for the volume functional.
\end{abstract}

\maketitle
 \markboth{Pengzi Miao  and Luen-Fai Tam}
 {Volume functional of compact manifolds with boundary}

\section{Introduction}
Given a compact $ n $-dimensional manifold $ \Omega$
with a boundary $ \Sigma$, we study variational
properties of the volume functional on the space
of constant  scalar curvature metrics on $ \Omega $ with
 a prescribed boundary metric on $ \Sigma$.
 The dimension $ n $ is assumed to be $  \geq 3 $.
 There are several
 motivations for us to consider this problem.

 The first motivation comes from a recent result  in
 \cite{FanShiTam07}.
 There one considers an asymptotically flat
 $ 3 $-manifold $(M, g)$ with a given end.
 Let $ \{ x_i \}$ be a coordinate system
 at $ \infty $ which defines the
 asymptotic structure of $ (M, g )$.
 Let  $ S_r = \{ x \in M \ | \ | x | = r \}  $
 be the coordinate  sphere, where $ | x | $
 denotes the coordinate  length.
  Let $ \gamma $
 be the induced metric on $ S_r $.
 When $ r $ is large, $ (S_r, \gamma)$ can be isometrically
  embedded in the Euclidean
 space $ \mathbb{R}^3 $ as a strictly convex hypersurface
 $ S^0_r$.  Let $ V_0 (r) $ be the volume of
  the region enclosed by $ S^0_r $ in $ \mathbb{R}^3 $
  and $ V_r $ be the volume of the region enclosed by $ S_r $
  in $ (M^3, g )$. It was proved in \cite{FanShiTam07} that,
 as $ r \rightarrow \infty $,
 \begin{equation} \label{sec0-eq-1}
 V ( r ) - V_0 (r) = 2 m_{ADM} \pi r^2 + o ( r^2 )
 \end{equation}
 whenever $ m_{ADM} $ is defined.  Here $ m_{ADM} $
 is the ADM mass of $ (M^3, g) $ \cite{ADM61}.
 Therefore if  $ R (g) $, the scalar curvature of $ g $,
  is nonnegative, then  by the Positive Mass Theorem
  \cite{SchoenYau79,Witten81}, one has
 \be \label{sec0-eq-2}
 V ( r ) \geq V_0 (r)
 \ee
 for sufficiently large $ r $. (Note that the case $ m_{ADM} = 0 $
 would imply $(M, g)$ is isometric to $ \mathbb{R}^3$, hence
 showing $ V(r) = V_0 (r) $.)

 As a statement on volume comparison, \eqref{sec0-eq-2} is
rather intriguing.  First, the curvature assumption it requires is
 only  on scalar curvature. Second, it was formulated as
 a boundary value problem, i.e. the competitors involved
 have a same Dirichlet boundary geometry.  It is natural
 to ask whether there exist related results on compact
 manifolds with boundary.

As a special case of Theorem \ref{BY-volume-s2-t3} and Theorem
\ref{global-minimum} in this paper, we have:

 \begin{thm} \label{sec0-thm-1}
 Let $ \Omega $ be a $ 3 $-dimensional compact manifold
 with a connected boundary $ \Sigma $. Let $ \gamma $ be a given
 metric on $ \Sigma $ such that $(\Sigma, \gamma)$ can be
 isometrically embedded in the Euclidean space $ \mathbb{R}^3 $
  as a  strictly convex hypersurface $ \Sigma_0 $.
   Let $ \mathcal{M}^0_\gamma $  be the
 space of zero scalar curvature metrics  $ g $ on $ \Omega $ such that
 the induced metric from $ g $ on $ \Sigma $ is $ \gamma $.
\begin{enumerate}
\item[(i)] Suppose
$ g \in \mathcal{M}^0_\gamma $ is a critical point of the
volume functional $ V(\cdot) $ on $ \mathcal{M}^0_\gamma$,
then
\begin{equation*}
V(g) \geq V_0,
\end{equation*}
where $ V_0 $ is the Euclidean volume of the region enclosed
by $ \Sigma_0 $ in $ \mathbb{R}^3$. Equality holds if and only
if $ (\Omega, g)$ is isometric to a standard Euclidean ball.
 \item[(ii)] There exists no element  in $ \mathcal{M}^0_\gamma$
 that  minimizes volume in  $ \mathcal{M}^0_\gamma$.
 \end{enumerate}
  \end{thm}

 Our second motivation comes from a variational characterization
 of Einstein metrics  on a closed manifold \cite{Besse87} \cite{Schoen89}:
Let $ \mathcal{M}_{-1} $ be the space of metrics with constant
scalar curvature $ - 1 $ on a compact manifold $ M $ without
boundary. Then an element $ g \in \mathcal{M}_{-1} $ is
 an Einstein metric if and only if $ g $ is a critical point of the
volume functional $ V(\cdot) $ on $ \mathcal{M}_{-1} $.
This characterization  follows from Proposition 4.47
in \cite{Besse87} or the argument preceding Lemma 1.2
in \cite{Schoen89}.
The proof used the fact that Einstein metrics correspond
to critical points of the total scalar curvature functional
$$ S (g) =   V(g)^\frac{2-n}{n}  \int_M R(g)  d V_g ,$$
while the restriction of $ S ( \cdot) $ to the space
$ \mathcal{M}_{-1} $ agrees with $ - V(\cdot)^\frac{2}{n} $.

We want to establish a similar result  on compact manifolds
with boundary, with an aim to therefore find a proper concept
of metrics that would sit between constant  scalar curvature
metrics and Einstein metrics.

As a special case of Theorem \ref{critical-nec-suff-s1t1} in this
paper, we have:
\begin{thm} \label{sec0-thm-2}
Let $ \Omega $ be a compact $ n $-dimensional manifold with smooth
boundary $ \Sigma $. Let $ \gamma $ be a given metric on $ \Sigma
$ and $ K= - n ( n - 1 )$ or $  0 $. Let $ \mathcal{M}^{K}_\gamma$
be the space of metrics on $ \Omega $ which have constant scalar
curvature $ K $ and have induced metric on $ \Sigma $ given by $
\gamma$. Then $ g \in \mathcal{M}^{ K }_\gamma $ is a critical
point of the volume functional $ V( \cdot ) $ on $
\mathcal{M}^K_\gamma$ if and only if there is a function $ \lambda
$ on $ \Omega $ such that $ \lambda = 0 $ on $ \Sigma $ and
\begin{equation}  \label{sec0-eq-3}
    -(\Delta_g\lambda)g+\nabla^2_g\lambda-\lambda \Ric(g)  = g,
\end{equation}
where $ \Delta_g $, $ \nabla^2_g $ are the Laplacian, Hessian
operator with respect to $ g $ and $ Ric(g) $ is the Ricci
curvature of $ g $.
\end{thm}

It will be shown in Theorem \ref{umbilic-s2-t2} that if $ g $ is a metric,
defined on some open set, which satisfies \eqref{sec0-eq-3}
with some function $ \lambda $, then $ g $ necessarily has
constant scalar curvature. Furthermore, if $ \Omega $ is
indeed a closed manifold and  $g $ is a metric on $ \Omega $
with negative scalar curvature  for
which  \eqref{sec0-eq-3} holds with some function
$ \lambda $, then $ g $ is an Einstein metric.

Our method in deriving the first variation of the volume functional
also leads to the following characterization of geodesic balls
in the hyperbolic  space $ \H^n $ and the sphere $ \mathbb{S}^n $
through  the first variation of the total mean curvature integral (see
Theorem \ref{critical-spaceforms-s2-t1} and Proposition
\ref{ball-in-HorS}).

\begin{thm} \label{sec0-thm-3}
Let $\Omega$ be a connected domain with compact closure
in $\mathbb{H}^n$ or $\mathbb{S}^n$ and with
a smooth (possibly disconnected)  boundary $\Sigma$.
If $\Omega \subset\mathbb{S}^n$, we also assume the volume
of $ \Omega $ is less than half of the volume of $ \mathbb{S}^n$. Let $g$
be the corresponding space form metric on $ \Omega $.
Then  $ \Omega $ is a geodesic ball if and only if
$$
\oint_\Sigma H^\prime (0)=0
$$ for any smooth variation $\{ g(t) \}$ of $g$ on $ \Omega$ such that
$ g(t) $ and $ g $ have the same scalar curvature and the same
induced boundary metric. Here $ H^\prime(0) $ is the variation of
the mean curvature of $ \Sigma $  in $(\Omega, g(t) )$ with
respect to the outward unit normal.
\end{thm}

In particular, Theorem \ref{sec0-thm-3} implies that if $ \Omega $ is
not a geodesic ball in $ \H^n $ or the half sphere $ \mathbb{S}^n_+$,
then there exists a deformation $ \tilde{g} $ of the standard
metric $ g $ on $ \Omega $ such that $ \tilde{g} $ and $ g $ have
the same scalar curvature and the same induced boundary metric
on $ \p \Omega $, and
$$ \oint_{\p \Omega} H < \oint_{\p \Omega} \tilde{H} ,$$
where $ H$, $ \tilde{H} $ are the mean curvature of $ \p \Omega $
in $(\Omega, g)$, $(\Omega, \tilde{g})$.
One wants to compare this with the results in \cite{ShiTam02}
and \cite{ShiTam07} (see the remark after Proposition \ref{ball-in-HorS}).

In \cite{Schoen89}, Schoen gave the following conjecture
concerning the volume functional $ V(\cdot)$ on a closed hyperbolic
manifold.

\vspace{.2cm}
\noindent {\bf Conjecture} {\em
Let $ (M^n, h )$ be a closed hyperbolic manifold. Let $ g $ be another
metric on $ M $ with $ R ( g) \geq R ( h ) $, then $ V(g) \geq V(h)$.}
\vspace{.2cm}

This conjecture remained widely open until recently its $ 3
$-dimensional case followed as a corollary of Perelman's work on
geometrization \cite{Perelman1,Perelman2}. It is natural to wonder
if there exists a similar  conjecture or result on compact
manifolds with boundary with a hyperbolic metric. In \cite{AST07},
related results were established by Agol, Storm and Thurston, 
on compact $ 3 $-manifolds whose
boundary are minimal surfaces. In this paper, we note that the
above conjecture on closed manifolds  does not generalize directly
to manifolds with boundary if only the Dirichlet boundary
condition is imposed.

 \begin{thm} \label{sec0-thm-4}
 Let $ \mathbb{H}^3 $ be the hyperbolic space and
 $ g_{\mathbb{H}} $ be the standard hyperbolic metric
 on $ \mathbb{H}^3$. There exists a small constant
 $ \delta > 0$ such that if $ B $ is a geodesic ball
 in $ \mathbb{H}^3 $ with geodesic radius less than $ \delta $,
 then there is another metric $ g $ on $ B $ such that
 $  g $ induces the same boundary metric on $ \p B $
 as $ g_{\mathbb{H} }$ does, $ R (g ) = R ( g_{\mathbb{H}}) $,
 and $ V( g ) < V( g_\mathbb{H} ) $ on $ B $.
 \end{thm}

 The paper is organized as follows. In Section \ref{sec-firstvar},
 we first analyze the manifold structure of the space of constant
 scalar curvature metrics with a prescribed boundary metric,
 then we compute the first variation of the volume functional
 and derive the critical point equation. We also relate the first
 variation of the volume functional on domains in an Einstein
 manifold with non-zero scalar curvature to the first variation
 of the total boundary mean curvature integral of the associated
 metric variation.
 In Section \ref{sec-spaceforms}, we show that the only
 domains in space forms, on which the standard metrics are
 critical points, are geodesic balls. We also establish some
 general  properties for metrics satisfying
 the critical point equation. Then we focus on the zero
 scalar curvature case to prove a theorem that
 compares the volume of a critical point with the
 corresponding Euclidean volume. In Section \ref{sec-secondvar},
 we derive the second variational formula and apply it to geodesic balls
 in space forms. In particular, we show that, on Euclidean balls
 and ``small'' hyperbolic and spherical balls in dimensions
 $ 3 \leq n \leq 5 $, the volume functional
 achieves a saddle point at the standard space form metrics.
 For completeness, we include an appendix in which we construct
 traceless and divergence free $ (0,2) $ symmetric tensors
 with prescribed compact support on space forms, which are needed
 in the second variation construction.\vskip .1cm

 {\it Acknowledgment}: The authors would like to thank
Robert Bartnik, Todd Oliynyk and Andrejs Treibergs
for useful discussions. 

\section{First variational formula for the volume  functional} \label{sec-firstvar}

Let $ \Omega $ be an $ n $-dimensional, connected,
compact manifold with a smooth (possibly disconnected)
boundary $\Sigma$.  Let $\gamma$ be a smooth
 metric on $\Sigma$. Let $\mathcal{S}^{k,2}$ be the space of
 $W^{k,2}$ symmetric (0,2) tensors on $\Omega$. We will
 always assume $k > \frac{n}{2} + 2$ so that
 each $h\in  \mathcal{S}^{k,2}$ is $C^{2,\alpha}$ up to the
 boundary. Let $\mathcal{S}^{k,2}_0$ be the subspace
 consisting those $h$ with $h|_{T(\Sigma)}=0$. That is to say,
 $h(v,w)=0$ for all $v, w$ tangent to $\Sigma$. Let
 $\mathcal{M}_{\gamma}$ be the open set in $ \mathcal{S}^{k,2}$
 consisting $g$ which is a Riemannian metric $g>0$ such that
 $g|_{T(\Sigma)}=\gamma$.
Let $W^{k-2,2}(\Omega)$ be the space of $W^{k-2,2}$ functions on
${\Omega}$. Let $\mathcal{R}$ be the scalar curvature map which
maps $g\in \mathcal{M}_{\gamma}$ to its scalar curvature $R(g)$. 
The fact that $ \mathcal{R} $ is a smooth map was shown by
Fischer and Marsden in \cite{FischerMarsden1975}:

\begin{lma}\label{s1-l1}
The map $\mathcal{R}:\mathcal{M}_{\gamma}\to W^{k-2,2}(\Omega)$ is
smooth.
\end{lma}

Let $g_0\in  \mathcal{M}_{\gamma}$. Suppose the scalar curvature
of $g_0$ is a constant $K$. We have

\begin{lma}\label{s1-l2}
Suppose $ 0 $ is not one of the Dirichlet eigenvalues of the operator
$$
(n-1)\Delta_{g_0}+K .
$$
 Then, near $g_0$, the set
$$
\mathcal{M}_\gamma^K=\{g\in \mathcal{M}_\gamma|\
\mathcal{R}(g)=K\}
$$
is a submanifold of $ \mathcal{M}_\gamma $.
\end{lma}
\begin{proof} The linearization of $\mathcal{R}$ at $g_0$
is
\begin{equation}\label{linear-scalar-e1}
D\mathcal{R}_{g_0}(h)=-\Delta_{g_0}(\tr_{g_0}h)+
\div_{g_0}\lf(\div_{g_0}(h)\ri)-\la h,\Ric(g_0)\ra_{g_0}
\end{equation}
for $h$ in the tangent space of $g_0$ in $\mathcal{M}_\gamma$,
 which is equal to $\mathcal{S}^{k,2}_0$. Here $\tr_{g_0} (h)$ and $ \div_{g_0} (h) $
 denote  the trace of $h$ and the divergence of $ h $ with respect to the metric $ g_0 $.
 Let $f \in W^{k-2,2}(\Omega)$ which is the tangent space at an element in
$W^{k-2,2}(\Omega)$. Consider the following boundary value
problem:
\begin{equation}\label{submanifold-e1}
    \begin{cases}-(n-1)\Delta_{g_0}u-Ku&=f\ \text{in
    $\Omega$}\\
    u&=0\ \text{on $\Sigma$},
    \end{cases}
\end{equation}
which has a unique solution by the assumption on the Dirichlet eigenvalues
and the Fredholm alternative, see
\cite{Fischer-ColbrieSchoen1980}. Let $h=ug_0$. We conclude that
$D\mathcal{R}_{g_0}(h)=f$. Hence $D\mathcal{R}_{g_0}$ is
surjective. As  $\mathcal{S}^{k,2}_0 $ is a Hilbert space, the
kernel of  $ D \mathcal{R}_{g_0} $ automatically splits. By the
implicit function theorem, the lemma follows.
\end{proof}

We want to apply Lemma \ref{s1-l2} to domains in space forms: the
Euclidean space $\mathbb{R}^n$, the hyperbolic space
$\mathbb{H}^n$ and the sphere $\mathbb{S}^n$.

\begin{cor}\label{s1-c1}
Suppose $\Omega$ is a connected domain with compact closure in
$\mathbb{R}^n$, $\mathbb{H}^n$ or $\mathbb{S}^n$ and with a
smooth (possibly disconnected) boundary $\Sigma$.
Let $g$ be the standard metric on $ \Omega $
and let $\gamma=g|_{\Sigma}$. Then, near $g$, the set
$\mathcal{M}_\gamma^K$ is a submanifold of
$\mathcal{M}_\gamma$, where (i) $K=0$ if $\Omega\subset
\mathbb{R}^n$, (ii) $K=-n(n-1)$ if $\Omega\subset
\mathbb{H}^n$, and (iii)  $K=n(n-1)$ if $\Omega\subset
\mathbb{S}^n$ and $V(\Omega)<\frac12 V(\mathbb{S}^n)$,
where $ V(\Omega) $ and $ V(\mathbb{S}^n)$ denote the
volume of $(\Omega, g)$ and $ \mathbb{S}^n$.
\end{cor}
\begin{proof} Cases (i) and (ii) follow immediately from Lemma
\ref{s1-l2}. In case (iii), since the volume of $\Omega\subset
\mathbb{S}^n$ is less than the
volume of a hemisphere, the first Dirichlet eigenvalue for the
Laplacian is larger than $n$ by the Faber-Krahn inequality
\cite{Sperner73}, the result follows from Lemma \ref{s1-l2}
again.
\end{proof}

Next we want to consider the volume functional
$$
V:\mathcal{M}_\gamma\to \mathbb{R},
$$
which is a smooth functional on $\mathcal{M}_\gamma$. For $g\in
\mathcal{M}_\gamma$, the first variation of $ V (\cdot) $ at $ g $ is:
\begin{equation}\label{1stvariation-e1}
  DV_g(h)=\frac12\int_{\Omega}\tr_g(h)dV_g .
\end{equation}
We are interested in critical points of $V(\cdot)$
restricted to $\mathcal{M}_\gamma^K $, when
 $\mathcal{M}_\gamma^K $ is a submanifold.
 First, we will construct explicit deformations on
 $\mathcal{M}_\gamma^K$.

\begin{prop}\label{deformation-p1}
Let $g_0\in \mathcal{M}_\gamma^K$ be a smooth metric such
that the first Dirichlet eigenvalue of
$(n-1)\Delta_{g_0}+K$ is positive. Let $h$ be a smooth
symmetric (0,2) tensor on $\Omega$  such that $h |_{ T (\p
\Omega)} = 0 $. Let $g(t)=g_0+th$ which is a smooth metric
provided $|t|$ is small enough. Then there is $t_0>0$ and
$\e>0$ such that for all $|t|<t_0$ the following Dirichlet
boundary value problem has a unique solution $u(t)$ such
that $1-\e\le u\le 1+\e$:
\begin{equation}\label{deformation-e1}
  \begin{cases}
    \alpha\Delta_{g(t)}u- {R}(t)u&=-Ku^a, \text{\ in $\Omega$}\\
    u&=1,  \text{\ on $\Sigma$,}
\end{cases}
    \end{equation}
where $\alpha=4(n-1)/(n-2)$, $a=(n+2)/(n-2)$, and ${R}(t)$
is the scalar curvature of $g(t)$. Moreover, $v=\frac{\p u}{\p
t}|_{t=0}$ exists and is a smooth function on $\ol{\Omega}$
(the closure of $ \Omega$), which is the unique solution of:
\begin{equation}\label{deformation-e2}
  \begin{cases}
    (n-1)\Delta_{g_0}v+Kv&= \frac{n-2}{4}
    \frac{\p {R}}{\p t}(0), \text{\ in $\Omega$}\\
    v&=0,  \text{\ on $\Sigma$.}
\end{cases}
    \end{equation}
\end{prop}
\begin{proof} Since the first eigenvalue of $(n-1)\Delta_{g_0}+K$
is positive, there is a smooth function $\phi$ on $\ol\Omega$ and
$\delta>0$ such that
\begin{equation}
    (n-1)\Delta_{g_0}\phi+K\phi+\delta \phi=0
\end{equation}
in $\Omega$ and $\phi>0$ on $\ol\Omega$, see
\cite{Fischer-ColbrieSchoen1980} for example.   Let
$ b =\min_{\ol\Omega}\phi$, then $ b > 0 $.
Suppose $1\ge t>0$ is small enough such that $g(t)$ is a
Riemannian metric on $\Omega$. For each $ t $,
let $ L(\cdot) $ denote the operator
$$ \alpha\Delta_{g(t)} (\cdot) - {R}(t) (\cdot)  + K (\cdot)^a . $$
Then
\begin{equation}\label{deformation-e3}
  \begin{split}
L(1+t\phi)&=\alpha \Delta_{g(t)}(1+t\phi)  - {R}(t)(1+t\phi)+K(1+t\phi)^a\\
&=\alpha\Delta_{g_0}(1+t\phi)-K(1+ t\phi)+
\alpha\lf(\Delta_{g(t)}-\Delta_{g_0}\ri)
\lf(1+t\phi\ri)\\
&\ \  -({R}(t)-K)(1+t\phi)+K(1+t\phi)^a\\
&\le -\frac{4}{n-2} \delta t\phi+C_1t+C_2t^2
\end{split}
\end{equation}
where $C_1$ is a positive constant depending only on $g_0$
and $h$, and $C_2>0$ depends also on $\phi$. Similarly,
\begin{equation}\label{deformation-e4}
  \begin{split}
L(1-t\phi)&=\alpha \Delta_{g(t)}(1-t\phi)  - {R}(t)(1-t\phi)+K(1-t\phi)^a\\
&=\alpha\Delta_{g_0}(1-t\phi)-K(1- t\phi)+
\alpha\lf(\Delta_{g(t)}-\Delta_{g_0}\ri)
\lf(1-t\phi\ri)\\
&\ \  -( {R}(t)-K)(1-t\phi)+K(1-t\phi)^a\\
&\ge \frac{4}{n-2} \delta t\phi-C_3t-C_4t^2
\end{split}
\end{equation}
where $C_3$ is a positive constant depending only on $g_0$
and $h$, and $C_4>0$ depends also on $\phi$. By rescaling
$\phi$, we may assume that $(4\delta b)/(n-2)\ge 2C_1$ and
$(4\delta b)/(n-2)\ge 2C_3$. Then for $t>0$ small enough
$L(1+t\phi)\le 0$ and $L(1-t\phi)\ge 0$. By \cite[Theorem
2.3.1]{Sa}, \eqref{deformation-e1} has a solution $u$
satisfying $1-t\phi\le u\le 1+t\phi$ provided $t>0$ is
small enough. The proof for $t<0$ is similar.

To prove uniqueness, let $u_1$ and $u_2$ be two solutions of
\eqref{deformation-e1} such that $1-\e\le u_1,u_2\le 1+\e$. Then
there is $\delta_1>0$ depending only on $g_0$ and $h$ such that
for $|t|$ small enough:
\begin{equation}\label{deformation-e5}
   \begin{split}
   \delta_1  \int_{\Omega}(u_1-u_2)^2dV_{g(t)}
   \le &  - \frac{4}{n-2}\int_\Omega K(u_1-u_2)^2dV_{g(t)}
   \\
    &  + \int_\Omega \lf[K(u_1^a-u_2^a) - {R}(t)(u_1-u_2)\ri](u_1-u_2)dV_{g(t)}\\
   &\le (C_5|t|+C_6\e)\int_\Omega (u_1-u_2)^2dV_{g(t)}
\end{split}\end{equation}
where $C_5$ and $C_6$ are constants depending only on $g_0$ and
$h$. Hence if $\e>0$ is small enough and $|t|$ is small enough, we
must have $u_1=u_2$.

To prove the last part of the proposition, let $u(t)$, $t\neq0$ be
the solution of \eqref{deformation-e1} obtained above,
then $ | u(t) - 1 | \leq | t | \phi $.
Let $w=(u-1)/t$. Then $w$ satisfies:
\begin{equation}\label{deformation-e6}
  \begin{cases}
    \alpha\Delta_{g(t)}w&= {R}(t)w-\frac{K - {R}(t)}t-
    \frac{ K(u^a-1) }t, \text{\ in $\Omega$}\\
      w&=0,  \text{\ on $\Sigma$.}
\end{cases}
    \end{equation}
    Since the right side of the equation is bounded by a constant
    independent of $t$ and $x\in \Omega$, there exist
    $\beta >0$ and $C_7>0$ independent of $t$ such that
    the H\"older norm with respect to $g_0$ with exponent $ \beta $
    of $w$ in $\Omega$ is bounded by $C_7$, see
    \cite[Theorem 8.29]{GilbargTrudinger}. By the Schauder
    estimates \cite[Theorem 6.19]{GilbargTrudinger}, one can
    conclude that for any $k\ge 2$, then $C^{k,\beta }$ norm of
    $w$ is uniformly bounded by a constant independent of $t$.
    Hence for any $t_i\to 0$, we can find a subsequence which
    converge to a solution $v$ of \eqref{deformation-e2}. Since
    $v$ is unique by the assumption on the first eigenvalue of
    $(n-1)\Delta_{g_0}+K$ and the Fredholm alternative, we
    conclude that $\frac{\p u}{\p t}|_{t=0}=v$. This completes the
    proof of the proposition.
\end{proof}

\begin{thm}\label{critical-nec-suff-s1t1}
Let $g \in \mathcal{M}_\gamma^K$ be a smooth metric such
that the first Dirichlet eigenvalue of
$(n-1)\Delta_{g }+K$ is positive.
Then $g$ is a critical point of the volume functional in
$\mathcal{M}_\gamma^K$ if and only if there is a smooth function
$\lambda$ on $\ol\Omega$ such that
\begin{equation}\label{critical-e1}
   \left\{%
\begin{array}{ll}
    -(\Delta_g\lambda)g+\nabla^2_g\lambda-\lambda \Ric(g)  &=g \hbox{\ in $\Omega$ ;} \\
    \lambda &=0 \hbox{\ on $\Sigma$.}
\end{array}%
\right.
\end{equation}
\end{thm}

\begin{proof} By Lemma \ref{s1-l2}, there is a neighborhood $U$ of
$g$ in $\mathcal{M}_\gamma$ such that
$U\cap \mathcal{M}^K_\gamma$ is a
submanifold. Suppose that $g$ is a critical point of $V(\cdot)$
in $\mathcal{M}_\gamma^K$. Since the first eigenvalue of
$(n-1)\Delta_g+K$ is positive, by the Fredholm alternative we have
a smooth function $\lambda$ on $\ol\Omega$ satisfying:
\begin{equation}\label{critical-e2}
\begin{cases}
    \Delta \lambda&=-\frac1{n-1}\lf(\lambda K+n\ri), \text{\ in $\Omega$}\\
    \lambda&=0,  \text{\ on $\Sigma$.}
\end{cases}
    \end{equation}
We want to prove that $\lambda$ satisfies the interior equation in
\eqref{critical-e1}.
 Let $h$ be a smooth symmetric (0,2) tensor with compact
support in  $\Omega$. Let $g(t)=g+th$, then $g (t) \in
\mathcal{M}_\gamma$ for small $ t $.
For each $t$, consider the following Dirichlet
boundary value problem:
\begin{equation}\label{critical-e3}
  \begin{cases}
    \alpha\Delta_{g(t)}u-R(t)u&=-Ku^a, \text{\ in $\Omega$}\\
    u&=1,  \text{\ on $\Sigma$,}
\end{cases}
    \end{equation}
where $\alpha=4(n-1)/(n-2)$, $a=(n+2)/(n-2)$ and $R(t)$ is the
scalar curvature of $g(t)$. For $|t|$ small, the equation has a
unique positive solution $u(t)$ which is smooth up to the boundary
by Proposition \ref{deformation-p1}. Moreover,
$u^{4/(n-2)}(t)g(t)$ is in $\mathcal{M}_\gamma^K$ and is a $C^1$
curve in $U\cap \mathcal{M}_\gamma^K$. Since $u\equiv 1$ at $t=0$,
$$
\frac{d}{dt}\lf(u^{4/(n-2)}(t)g(t)\ri)|_{t=0}=\frac4{n-2} u'g+h,
$$
where $u'=u'(0)$. Since $g$ is a
critical point of $V(\cdot)$ in $ \mathcal{M}_\gamma^K$,
by \eqref{1stvariation-e1}, we have
\begin{equation}\label{critical-e4}
    \int_\Omega \lf(\frac{4n}{n-2}u'+\tr_g(h)\ri)dV_g=0.
\end{equation}
 Now by Proposition \ref{deformation-p1} again, $u'$ satisfies:
\begin{equation}\label{critical-e5}
  \begin{cases}
    \alpha\Delta_{g}u'-R'(0)-Ku'&=-aKu', \text{\ in $\Omega$}\\
    u'&=0,  \text{\ on $\Sigma$.}
\end{cases}
    \end{equation}
Since $\lambda=0$ on $\Sigma$, by \eqref{critical-e2} and
\eqref{critical-e5}, we have
\begin{equation}\label{critical-e6}
    \begin{split}
  \frac{4n}{n-2}\int_\Omega u' dV_g&= \frac{n\alpha}{n-1}\int_\Omega u' dV_g \\
        &= \alpha\int_\Omega \lf(-\frac{\lambda Ku'}{n-1}-u' \Delta_g
        \lambda\ri)dV_g\\
        &=\alpha\int_\Omega \lf(-\frac{\lambda Ku'}{n-1}-\lambda \Delta_g
        u'\ri)dV_g\\
        &=\int_\Omega \lf(-\frac{\alpha\lambda Ku'}{n-1}-\lambda R'(0)-(1-a)K\lambda u'\ri)dV_g\\
        &= \int_\Omega \lambda\lf(\Delta_{g}(\tr_{g}h)-
\div_{g}\lf(\div_{g}(h)\ri)+\la h,\Ric(g)\ra_{g}\ri)dV_g
\end{split}
\end{equation}
where we have used \eqref{linear-scalar-e1} in the last step.
Combining this with \eqref{critical-e4} we have
\begin{equation}\label{critical-e7}
  \int_\Omega\lambda \lf( \Delta_{g}(\tr_{g}h)-
\div_{g}\lf(\div_{g}(h)\ri)+\la h,\Ric(g)\ra_{g} \ri) +\tr_g(h) dV_g=0.
\end{equation}
Let $f=\tr_g(h)$. Let $\nu$ be the unit outward normal of
$\Sigma$. Since $f$ has compact support in $\Omega$, we have
\begin{equation}\label{critical-e8}
\int_\Omega\lambda\Delta_gfdV_g=\int_\Omega f\Delta_{g}\lambda dV_g
+\oint_\Sigma(\lambda f_\nu- f\lambda_\nu)=
\int_\Omega f\Delta_{g}\lambda dV_g
\end{equation}
where we use $ \psi_\nu $ to denote $ \frac{\p \psi}{\p \nu} $ for
a smooth function $ \psi $ on $ \bar{\Omega} $, and
\begin{equation}\label{critical-e9}
\begin{split}
    \int_\Omega \lambda \div_{g}\lf(\div_{g}(h)\ri) &=\int_\Omega \lambda g^{ij}g^{kl}h_{ik;jl}dV_g\\
&=- \int_\Omega \lambda_{l} g^{ij}g^{kl}h_{ik;j}dV_g+ \oint_\Sigma\lambda (\div h)_k\nu^k\\
&=\int_\Omega \la\nabla^2_g\lambda,h\ra dV_g
- \oint_\Sigma h(\nabla \lambda,\nu)+\oint_\Sigma\lambda (\div h)_k\nu^k\\
&=\int_\Omega \la\nabla^2_g\lambda,h\ra dV_g
\end{split}
 \end{equation}
 where we have used the fact that $h$ has compact support.
 Combining \eqref{critical-e7}--\eqref{critical-e9}, we have:
\begin{equation}\label{critical-e10}
0 =\int_\Omega  \la h,  \lf(\Delta_{g}\lambda\ri)g-\nabla^2_g\lambda+\lambda\Ric(g)+
   g\ra_g dV_g.
 \end{equation}
Since $h$ is arbitrary, $\lambda$ must satisfy the interior equation in
\eqref{critical-e1}.

To prove sufficiency, suppose there is a smooth function $\lambda$
satisfying \eqref{critical-e1}, let $h\in  \mathcal{S}^2_0$ be
in the tangent space of $g$ in $\mathcal{M}_\gamma^K$. Then $h$ is
in the kernel of $D\mathcal{R}_{g}$. Let $f=\tr_g(h)$ as before.
By \eqref{linear-scalar-e1} and the computation in
\eqref{critical-e8}--\eqref{critical-e9}, we have
\begin{equation}\label{critical-e11}
    \begin{split}
0&=\int_\Omega\lambda\lf(\Delta_g f-\div_g\lf(\div_g(h)\ri)+\la h,\Ric(g)\ra_g\ri)dV_g\\
&=\int_\Omega  \la h,  \lf(\Delta_{g}\lambda\ri)g-\nabla^2_g\lambda+\lambda\Ric(g)+
   g\ra_g dV_g - \int_\Omega fdV_g\\
   &\ \ -\oint_\Sigma f\lambda_\nu + \oint_\Sigma h(\nabla \lambda,\nu)
   + \oint_\Sigma \lambda \lf[ f_\nu - \lf( \div_g h \ri)(\nu) \ri] \\
   &=-\int_\Omega fdV_g\\
   &=-2DV(h),
\end{split}
\end{equation}
where we have used the fact that $\lambda$ satisfies
\eqref{critical-e1}, $h|_{T(\Sigma)}=0$, and $\lambda=0$ on
$\Sigma$. Hence $g$ is a critical point of $V(\cdot)$ in
$ \mathcal{M}_\gamma^K $.
\end{proof}

\begin{rmk}\label{critical-nec-suff-remark}
\begin{enumerate}
\item[(i)] From the proof of Theorem \ref{critical-nec-suff-s1t1},
    one can see that
    under the assumptions on $g$, $g$ is a critical point of
    the volume functional $ V (\cdot) $ in $\mathcal{M}_\gamma^K$ if and
    only if $V'(0)=\frac{d}{dt}V(g(t))|_{t=0}=0$ for any
    smooth variation $\{ g(t) \}$ of $ g $ in $\mathcal{M}_\gamma^K$.
 \item[(ii)]     The differential equation in
     \eqref{critical-e1} can be equivalently written as $ D
     \mathcal{R}^*_g ( \lambda ) = g $, where $ D
     \mathcal{R}^*_g $ is the formal $ L^2 $-adjoint of $ D
     \mathcal{R}_g $, and  a weak form of \eqref{critical-e1}
     can also be derived using the infinite dimensional
     Lagrangian multiplier method  employed by Bartnik  in
     \cite{Bartnik2005}.
\end{enumerate}
\end{rmk}

 Theorem \ref{critical-nec-suff-s1t1} shows that, for a constant
 scalar curvature metric $ g $  to be a critical point for $ V(\cdot)$
 in  $ \mathcal{M}^K_\gamma $, there need to exist  a function
 $\l $ which satisfies both the interior equation
\begin{equation} \label{interior-eq1}
  -(\Delta_g\lambda)g+\nabla^2_g\lambda-\lambda \Ric(g)  =g \ \
  \mathrm{on} \ \Omega
\end{equation}
and the boundary condition $ \l |_{\Sigma} = 0 $.
(Later in Theorem \ref{umbilic-s2-t2} one will see that
(\ref{interior-eq1})  alone  implies that $ g $ has constant scalar
 curvature.) It remains interesting to know what the first variation
 of $ V(\cdot)$ in $ \mathcal{M}^K_\gamma $ would be
 if only (\ref{interior-eq1}) is satisfied but $ \l |_{\Sigma} $ is not
 necessarily zero.

\begin{prop} \label{lnonzero}
Let $g \in \mathcal{M}_\gamma^K$ be a smooth metric. Suppose
there exists a smooth function $ \l $ on $ \bar{\Omega} $ such that
\begin{equation} \label{interior-eq2}
  -(\Delta_g\lambda)g+\nabla^2_g\lambda-\lambda \Ric(g)  =g \ \
  \mathrm{on} \ \Omega.
\end{equation}
Let $ \{ g(t) \} $ be a smooth path of metrics in $
\mathcal{M}_\gamma $ such that $ g(0) = g $ and $ g(t) \in
\mathcal{M}^K_\gamma $. Then
\begin{equation}
\frac{ d}{dt} V( g(t) ) |_{ t = 0 } = \oint_\Sigma \l  H^\prime (0) ,
\end{equation}
where $ H = H (t) $ is the mean curvature of $ \S $ in $ (\Omega,
g(t) )$ with respect to the unit outward pointing normal vector $
\nu $.
\end{prop}

\begin{proof}
Similar to (\ref{critical-e11}), we have
\begin{equation}
    \begin{split}
0 &=\int_\Omega\lambda\lf(\Delta_g f-\div_g\lf(\div_g(h)\ri)+\la h,\Ric(g)\ra_g\ri)dV_g\\
&=\int_\Omega  \la h,  \lf(\Delta_{g}\lambda\ri)g-\nabla^2_g\lambda+\lambda\Ric(g)+
   g\ra_g dV_g - \int_\Omega fdV_g\\
   &\ \ -\oint_\Sigma f\lambda_\nu + \oint_\Sigma h(\nabla \lambda,\nu)
   + \oint_\Sigma \lambda \lf[ f_\nu - \lf( \div_g h \ri)(\nu) \ri] ,
\end{split}
\end{equation}
where we let $ h = g^\prime ( 0 ) $, $ f = \tr_g h $ and used the fact
$ D \mathcal{R}_{g} ( h ) = 0 $.
By (\ref{interior-eq2}), we have
\begin{equation} \label{critical-e12}
    \begin{split}
\int_\Omega fdV_g =
 -\oint_\Sigma f\lambda_\nu + \oint_\Sigma h(\nabla \lambda,\nu)
   + \oint_\Sigma \lambda \lf[ f_\nu - \lf( \div_g h \ri)(\nu) \ri] .
\end{split}
\end{equation}
Let $ \nabla_\Sigma \l $ be the gradient of $ \l $ on $ (\Sigma, \gamma)$
and integrate by parts
\begin{equation}
\begin{split} \label{critical-e13}
 \oint_\Sigma h(\nabla \lambda,\nu) = & \oint_\Sigma h(\nabla_\Sigma \lambda,\nu)
 + \oint_\Sigma h(\nu, \nu) \l_\nu \\
 = &  - \oint_\Sigma \l \div_\gamma X
 + \oint_\Sigma h(\nu, \nu) \l_\nu,
\end{split}
\end{equation}
where $ X $ is the vector field on $ \Sigma $ that is dual
to the one form $ h(\nu, \cdot)|_{T (\Sigma)}$ on $
(\Sigma, \gamma)$ and $ \div_\gamma X $ denotes the
divergence of $ X $ on $ (\Sigma, \gamma)$. Plug
(\ref{critical-e13}) in (\ref{critical-e12}), we have
\begin{equation} \label{critical-e14}
    \begin{split}
\int_\Omega fdV_g =
 -\oint_\Sigma \la \gamma, h \ra_\gamma   \lambda_\nu
   + \oint_\Sigma \lambda \lf[ - \div_\gamma X +  f_\nu - \lf( \div_g h \ri)(\nu) \ri] .
\end{split}
\end{equation}
Now, let $ p \in \Sigma $ and let $ \{ x^i \ | \ i = 1, \ldots, n \}$ be
a coordinate chart around $ p $ in $ \Omega $ such that
$\{ x^A \ | \ A = 1 , \ldots, n-1 \}$ gives a  coordinate chart on $ \Sigma $
and $ \partial_n = \nu $.  Direct calculation shows
\begin{equation*}
( \div_g h )_n = h_{nn; n} + \gamma^{AB} h_{ A n ;B},
\end{equation*}
\begin{equation*}
\begin{split}
h_{nA;B} = & \p_B ( h_{nA} ) - \Gamma^C_{AB} h_{nC} - \Gamma^n_{AB} h_{nn}
- \Gamma^i_{nB} h_{iA} \\
= & X_{A;B} + \Pi_{AB} h_{nn} - \Pi^C_B h_{CA},
\end{split}
\end{equation*}
where $ \Pi_{AB} = \la \nabla_{\p_A} \nu , \p_B \ra_\gamma $ is the second
fundamental form of $ \Sigma $. Thus,
\begin{equation} \label{divXanddivh}
( \div_g h )_n = h_{nn; n} + \div_\gamma X + H h_{nn} - \la \Pi, h \ra_\gamma .
\end{equation}
 On the other hand,
we have the following formula of the linearization
of the mean curvature, see equation (42) in \cite{static03} for example:
\begin{equation} \label{derivativeH}
H^\prime (0) =  \frac12 h_{nn;n} + \frac12 H h_{nn} - \la \Pi, h \ra_g
- [ \div_g h - \frac12 d ( \tr_g h ) ]_n .
\end{equation}
(Note that the sign convention of $\Pi$ in \cite{static03} is
opposite to the one used here.)
Hence (\ref{divXanddivh}) and (\ref{derivativeH})  imply that
\be
2 H^\prime (0 ) = [  d ( \tr_g h )  - \div_g h ]_n - \div_\gamma X
- \la \Pi, h \ra_\gamma .
\ee
Therefore, (\ref{critical-e14}) becomes
\begin{equation} \label{critical-e15}
    \begin{split}
\int_\Omega fdV_g =
 -\oint_\Sigma \la \gamma, h \ra_\gamma   \lambda_\nu
   + \oint_\Sigma \lambda \lf[  2 H^\prime(0) +  \la \Pi, h \ra_\gamma  \ri] .
\end{split}
\end{equation}
Finally, by  \eqref{linear-scalar-e1} and the boundary condition $ h|_{T (\Sigma)} = 0 $,
we have
\be
\frac{d}{dt} V (g(t)) |_{t=0} =  \oint_\Sigma \l H^\prime (0) .
\ee
\end{proof}

As an application, we have the following:

\begin{cor}\label{lnonzero-cor}
Suppose $g \in \mathcal{M}_\gamma^K$ is a smooth Einstein metric
with $\Ric (g) = \kappa g$ (so $K=n\kappa$).
Suppose $ \kappa \neq 0$.
Let $ \{ g(t) \} $ be a smooth path of metrics in $
\mathcal{M}_\gamma $ such that $ g(0) = g $ and $ g(t) \in
\mathcal{M}^K_\gamma $. Then
\begin{equation}\label{lnonzero-cor-e1}
\frac{d}{dt} V(g(t)) |_{t=0} =-\frac1\kappa \oint_\Sigma H'(0) .
\end{equation}
In particular, if the first Dirichlet eigenvalue of
$(n-1)\Delta_g+K$ is positive, then $g$ is a critical point of the
volume functional $ V ( \cdot ) $ in $ \mathcal{M}^K_\gamma$
if and only if $ \oint_\Sigma H'(0)=0 $ for   any
smooth variation $ \{ g(t) \} $ of $g$   in $ \mathcal{M}_\gamma^K$.
\end{cor}

\begin{proof} If $g$ satisfies $\Ric (g) =\kappa g$, then the constant
function $\lambda=-1/\kappa$ satisfies \eqref{interior-eq2}. By
Proposition \ref{lnonzero}, \eqref{lnonzero-cor-e1} is true. The
last part of the corollary follows from \eqref{lnonzero-cor-e1}
and Theorem \ref{critical-nec-suff-s1t1}, see Remark
\ref{critical-nec-suff-remark}(i).
\end{proof}

For example, if $\Omega$ is a domain in $\mathbb{S}^n$, then
$$
\frac{d}{dt} V (g(t)) |_{t=0} = -\frac1{n-1} \oint_\Sigma H'(0)
$$
for any smooth variation $ \{ g(t) \} $ of the standard metric on $\Omega$
which keeps the induced boundary metric  fixed.

 \section{Critical points in space forms} \label{sec-spaceforms}
In this section, we shall discuss  the volume
 functional on domains in space forms.

 \begin{thm} \label{critical-spaceforms-s2-t1}
Let $\Omega$ be a connected domain with compact closure
in $\mathbb{R}^n$, $\mathbb{H}^n$ or $\mathbb{S}^n$ and with
a smooth (possibly disconnected)  boundary $\Sigma$.
If $\Omega\subset\mathbb{S}^n$, we also assume that
$V(\Omega)<\frac12V(\mathbb{S}^n)$. Let $g$ be the standard metric
on $ \Omega $ and let $\gamma = g |_{ T(\Sigma) }$. Suppose $g$ is a critical
point of the volume functional $ V (\cdot) $ in $\mathcal{M}_\gamma^K$, where
(i) $K=0$ if $\Omega\subset \mathbb{R}^n$, (ii) $K=-n(n-1)$ if
$\Omega\subset \mathbb{H}^n$, and (iii)  $K=n(n-1)$ if
$\Omega\subset \mathbb{S}^n$, then $ \Omega $ is a geodesic
ball. Conversely, if $\Omega$ is a geodesic ball, then the
standard metric $ g $ is a critical point of the volume functional
$ V ( \cdot ) $ in $\mathcal{M}_\gamma^K$.
 \end{thm}

 \begin{proof}
 By Theorem \ref{critical-nec-suff-s1t1},
 $(\Omega,g)$ is a critical point if and only if there exists a smooth function
 $\lambda$  such that
 \begin{equation*}
   \left\{%
\begin{array}{ll}
    -(\Delta_g\lambda)g+\nabla^2_g\lambda-\lambda \Ric(g)  &= g \hbox{\ in $\Omega$ ;} \\
    \lambda &=0 \hbox{\ on $\Sigma$.}
\end{array}%
\right.
\end{equation*}
Taking trace of the equation, we have:
\begin{equation}\label{critical-s2-e1}
    \Delta_g \lambda=-\frac1{n-1}\lf(\lambda K+n\ri)=\left\{%
\begin{array}{ll}
     -\frac{n}{n-1}, & \hbox{\ if $K=0$;} \\
      \lambda n  -\frac{n}{n-1},&   \hbox{\ if $K=-n(n-1)$;} \\
      -\lambda n  -\frac{n}{n-1}, &\hbox{\ if $K=n(n-1)$.} \\
\end{array}%
\right.
\end{equation}
Hence $\lambda$ satisfies
\begin{equation}\label{critical-s2-e2}
    \nabla^2_g \lambda=\lambda \Ric(g)-\frac{ \lambda K+1 }{n-1} g=\left\{%
\begin{array}{ll}
     -\frac{1}{n-1}g, & \hbox{\ if $K=0$;} \\
     \lf( \lambda    -\frac{1}{n-1}\ri)g,&   \hbox{\ if $K=-n(n-1)$;} \\
      \lf(-\lambda    -\frac{1}{n-1}\ri)g, &\hbox{\ if $K=n(n-1)$.} \\
\end{array}%
\right.
\end{equation}

In case (i), where $ \Omega $ is a domain in $ \mathbb{R}^n$,
(\ref{critical-s2-e2}) directly implies
\begin{equation} \label{linR}
 \lambda=-\frac 1{2(n-1)}|x|^2 +\sum_{i=1}^nb_ix^i+c, 
\end{equation}
where $x^1,\cdots, x^n$  are the standard coordinates on
$ \mathbb{R}^n$ and $b_i$, $c$ are constants.
By translating the origin, we may assume
that $\lambda=-\frac 1{2(n-1)}|x|^2+c$ for a possibly different
$c$. Since $\lambda$ is zero at the boundary of $\Omega$, $\Omega$
must be a Euclidean ball  and $c=\frac 1{2(n-1)}R^2$, where $R$ is
the radius of the ball. Conversely, if $\Omega$ is a Euclidean
ball of radius $R$, then it is easy to see that $\lambda$ given
above satisfies the conditions in Theorem \ref{critical-nec-suff-s1t1},
therefore the standard metric is a critical point.

Next we consider case (ii), where $\Omega$ is a domain in $\H^n$.
Suppose the standard metric is a critical point of the volume
functional. Let $\lambda$ satisfy  \eqref{critical-s2-e2},
then $\lambda$ is not identically zero. Since $\lambda=0$
on $ \Sigma$, there must be an interior point $ p \in \Omega$
such that $\nabla \lambda(p)=0$. Henceforth, we use
$ \nabla  $ to denote the covariant derivative with respect to $ g $.

 Embed $\mathbb{H}^n$ in $\mathbb{R}^{n,1}$, the Minkowski
 space with metric $dx_1^2+\cdots+dx_n^2-dt^2$ such that
$$
\mathbb{H}^n=\{(x_1,\dots,x_n,t) \ | \ x_1^2+\cdots+x_n^2-t^2=-1, \ t > 0 \}
$$
and such that $p$ is mapped to the point $(0,\cdots,0,1)$. Then
$\nu=\sum_{i=1}^n x_i \frac{\p}{\p x^i}+t\frac{\p}{\p t}$ is normal
to $\mathbb{H}^n$ with $\la \nu,\nu\ra=-1$ in $\mathbb{R}^{n,1}$.
Let $\nabla$ and $D$ be the covariant derivatives of $\mathbb{H}^n$
and $\mathbb{R}^{n,1}$ respectively. Consider a point in
$\mathbb{H}^n$ and a function $f$ defined near that point in
$\mathbb{R}^{n,1}$. Let $\{ e_i \} $ be a basis for the tangent space of
$\H^n$, we have
$$
\nabla_i\nabla_j f=D_iD_jf+ \Pi_{ij}  \la\grad f, \nu\ra
$$
where $\grad f$ is the gradient of $f$ in $\R^{n,1}$, $\Pi_{ij}=-\la
D_{e_i}e_j, \nu \ra$ is the second fundamental form of $\H^n$ and is
equal to the induced metric $g_{ij}$ on $\H^n$. Let $f$ be the
function $at+\frac1{n-1}$ where $a$ is chosen so that
$a+\frac1{n-1}=\lambda(p)$. Hence $f(p)=\lambda(p)$ and $\nabla
f(p)=0$. Since $\la \grad f ,\nu\ra=at$, one can check that $f$
satisfies \eqref{critical-s2-e2}. So $\nabla^2(\lambda-f)=(\lambda-f)g$.
Consider a geodesic $ \sigma (s) $ on $\H^n$ emanating from
$p$. Restricted to $ \sigma(s) $,
the function $\lambda-f$ satisfies the ODE
$$
(\lambda-f)''= \lambda-f,
$$
where `` $ ' $ " is taken with respect to $ s $. Since initially
$\lambda-f=(\lambda-f)'=0$, $ \lambda-f $ must be identically
zero. So $\lambda=f$. But $f=\lambda=0$ at the boundary of
$\Omega$. Hence $\Omega$ is a geodesic ball. When restricted on
$\H^n$, $t=\cosh r$ where $r$ is the geodesic distance from the
point $(0,\cdots,0,1)$. Hence $a=-\lf((n-1)\cosh R\ri)^{-1}$ where
$R$ is the radius of the geodesic ball $\Omega$.

Conversely, if $\Omega$ is a geodesic ball with center at
$(0,\cdots,0,1)$ with geodesic radius $R$. Then 
\begin{equation} \label{linH}
\lambda= \frac1{n-1}\lf(1-\frac{\cosh r}{\cosh R} \ri)
\end{equation}
 satisfies the conditions in Theorem \ref{critical-nec-suff-s1t1}, which implies the
standard metric is a critical point.

Finally we consider case (iii), where $\Omega$ is a domain with compact
closure in $ \mathbb{S}^n$. Suppose the standard metric is a
critical point of the volume functional. Let $\lambda$ satisfy
\eqref{critical-s2-e2}, then $\lambda$ is not identically zero.
 Since $\lambda=0$ on $ \Sigma $, there must be a
 point $ p \in \Omega $ such that $ \nabla \lambda (p) = 0 $.
 Embed $\mathbb{S}^n$ in
 $\mathbb{R}^{n+1} = \{ (x_0, x_1, \ldots, x_n ) \}$
 as the unit sphere  centered at the origin such that
 $p$ is mapped to the point $(0,\cdots,0,1)$.
 Let $f = a x_n  - \frac1{n-1}$, where $a$ is chosen so that $a - \frac1{n-1}
 =\lambda(p) $.
Hence $f(p)=\lambda(p)$ and $\nabla f(p)=0$. As in case (ii), one
can check that $f$ satisfies \eqref{critical-s2-e2}. Then one can
prove as before that $\lambda=f$. Since $ f = \lambda = 0 $
on the boundary of $ \Omega$, $ \Omega $
must be a geodesic ball centered at $ (0, \ldots, 0, 1)$
whose boundary can not be the equator
$ \{ x_n = 0 \} $. As we assume $V(\Omega)<\frac12V(\mathbb{S}^n)$,
$ \Omega $ must be strictly contained in the upper hemisphere.
When restricted to $ \mathbb{S}^n$, $ x_n = \cos{ r }$
where $r$ is the geodesic
distance from the point $(0,\cdots,0,1)$. Hence $a= \lf((n-1)\cos
R \ri)^{-1}$ where $R$ is the radius of the geodesic ball $\Omega$.

  Conversely, if $\Omega$ is a geodesic ball in $ \mathbb{S}^n $
  with radius $R \neq \frac{\pi}{2}$, then
  \begin{equation} \label{linS}
   \lambda = \frac1{n-1}\lf( \frac{\cos{r} }{\cos{R} } - 1 \ri)
  \end{equation} 
   satisfies \eqref{critical-e1} in Theorem \ref{critical-nec-suff-s1t1}.
  If furthermore $ \Omega $ is contained in a hemisphere, then
  Theorem \ref{critical-nec-suff-s1t1} shows that
  the standard metric is a critical point.
\end{proof}

As a direct application of Theorem \ref{critical-spaceforms-s2-t1} and Corollary
\ref{lnonzero-cor}, we have

\begin{prop} \label{ball-in-HorS}
Let $\Omega$ be a connected domain with compact closure
in $\mathbb{H}^n$ or $\mathbb{S}^n$ and with
a smooth (possibly disconnected)  boundary $\Sigma$.
If $\Omega \subset\mathbb{S}^n$, we also assume that
$V(\Omega)<\frac12V(\mathbb{S}^n)$. Let $g$ be the standard
metric on $ \Omega $, let $ K $ be the constant that is equal to
the scalar curvature of $ g $ and let $\gamma = g |_{ T(\Sigma) }$.
Then  $ \Omega $ is a geodesic ball if and only if
$$
\oint_\Sigma H^\prime (0)=0
$$ for any smooth
  variation $\{ g(t) \}$ of $g$ in $ \mathcal{M}_\gamma^K$. Here
  $ H^\prime(0) $ is the variation of the mean curvature of $ \Sigma $
  in $(\Omega, g(t) )$ with respect to the outward unit normal.
  \end{prop}

\begin{proof} By Theorem \ref{critical-spaceforms-s2-t1}, $ \Omega $
is a geodesic ball if and only if the standard metric $ g $ is a critical point
of the volume functional $ V (\cdot) $ in $ \mathcal{M}^K_\gamma$. On the other
hand, as $ g $ is an Einstein metric with non-zero scalar curvature,
 Corollary \ref{lnonzero-cor} shows $ g $ is a critical point
of $ V(\cdot) $ in $ \mathcal{M}^K_\gamma$ if and only if
$$
\oint_\Sigma H^\prime (0)=0
$$ for any smooth variation $\{ g(t) \}$ of $g$ in $ \mathcal{M}_\gamma^K$.
The proposition follows.
 \end{proof}

Before we proceed to discuss properties of general critical metrics of
$ V(\cdot)$, we want to relate the result in Proposition
\ref{ball-in-HorS} to the results in \cite{ShiTam02} and \cite{ShiTam07}.
Let $ \Sigma $ be any given compact strictly convex hypersurface in
$ \R^3$. If $(\Omega,g)$ is a compact Riemannian $ 3$-manifold with
nonnegative scalar curvature whose boundary is isometric to $ \Sigma $
and has positive mean curvature $ H$, then it was proved
in \cite{ShiTam02} that
\begin{equation}\label{H-e1}
\oint_\Sigma H_0 \ge \oint_{\p \Omega} H,
\end{equation}
where $H_0$ is the mean curvature of $\Sigma$ in $\R^3$.
In particular, if $ \Omega $ is the domain enclosed by
$ \Sigma $ in $ \R^3$,   \eqref{H-e1} then implies that
 $$ \int_\Sigma H'(0)=0 $$
  for any smooth variation $ \{g(t) \}$ of $ g $ in $ \mathcal{M}^0_\gamma$.
This contrasts sharply with Proposition \ref{ball-in-HorS}, by which
we know the unique compact convex surfaces in $ \H^3 $ or $ \mathbb{S}^3_+$
with that property are geodesic spheres.
As a result, Proposition \ref{ball-in-HorS} implies that \eqref{H-e1}
does not generalize directly to an arbitrary compact convex surface
in $ \H^3 $ or $ \S^3_+$.  On the other hand, it was proved
in \cite{ShiTam07} that  \eqref{H-e1} does hold if $ \Sigma $ is a
geodesic sphere in $ \H^3$ and $ (\Omega, g)$ has scalar curvature
no less than $ -6 $, which is consistent with Proposition \ref{ball-in-HorS}.
It remains an interesting question to know whether \eqref{H-e1} is true for
a compact $ 3$-manifold $(\Omega, g)$ whose boundary
is isometric to a geodesic sphere in $ \S^3_+$ and whose
scalar curvature is greater than or equal to $ +6$.

Next we discuss some general properties of a Riemannian metric $ g $
for which there exists a function $\lambda$ satisfying the differential equation
in \eqref{critical-e1} in Theorem \ref{critical-nec-suff-s1t1}.

\begin{thm}\label{umbilic-s2-t2} Let $(\Omega,g)$ be a
connected, smooth Riemannian manifold. Suppose there is a
smooth function $\lambda$ on $\Omega$ such that
\begin{equation}\label{umbilic-e1}
    -(\Delta_g \lambda)g+\nabla^2_g \lambda -\lambda \Ric(g) =g.
\end{equation}
Then
\begin{enumerate}
\item[(i)] $g$ has constant scalar curvature.
\vspace{.1cm}
\item[(ii)] If $ \Omega $ is compact without boundary and
$ g $ has negative scalar curvature, then $ g $ is an Einstein metric.
\vspace{.1cm}
\item[(iii)] If $\Omega$ is compact with a smooth (possibly disconnected)
         boundary $\Sigma$ such that $\lambda = 0 $ at $ \Sigma$,
         and if
         the first Dirichlet eigenvalue of $(n-1)\Delta_g +K$ is nonnegative
         where $K$ is the scalar curvature of $g$, then
         along each connected component $ \Sigma_\alpha $ of $ \Sigma$,
         the mean curvature of $ \Sigma_\alpha $
         with respect to the
         outward unit normal $ \nu $ is a positive constant. In fact,
         $ \Sigma_\alpha $ is umbilic and  its second fundamental form $ \Pi_\alpha $
         satisfies $\Pi_\alpha = a_\alpha g |_{ T (\Sigma_\alpha ) }$ for some
         constant $ a_\alpha > 0 $.
\item[(iv)]  Under the same assumptions as in (iii), at each point in $ \Sigma$,
          \begin{equation}\label{ric-boundary-e}
             2\Ric(\nu,\nu)+R_\Sigma=K+\frac{n-2}{n-1}H^2 ,
         \end{equation}
where  $R_\Sigma$  is the scalar curvature of $\Sigma$ and $ H $ is
mean curvature of $ \Sigma $.
\end{enumerate}
\end{thm}

\begin{proof} (i) can be proved as in \cite{FischerMarsden1975}.  By
taking the divergence of \eqref{umbilic-e1}, and using the Bianchi
identity, we conclude that $\lambda d R =0$ where
${R}$ is the scalar curvature of $g$. At a point $p$ where
$\lambda(p)\neq0$, we have $d {R} =0$. Suppose $ p \in
\Omega $ is an interior point where $ \l (p ) = 0 $. The equation
 \begin{equation} \label{eq-trace}
\Delta_g \l = - \frac{n}{n-1} - \frac{{R}}{n-1} \l
\end{equation}
 then implies $ \triangle \l ( p ) < 0 $. Thus, either $ \nabla \l (p) \neq 0 $ or
$ p $ is a  strict local maximal point for $ \l $. In either case,
we would have $ d {R} = 0 $ in a neighborhood around $p$.
Hence, $d R =0 $ in $ \Omega$.

To prove (ii), we know that $ R $ is a negative constant
by (i) and the assumption. As $ \Omega $ is compact
without boundary, it follows from \eqref{eq-trace}
and the maximum principle that $ \lambda $ must
be a constant. Hence, $ g $ is an Einstein metric.

To prove (iii), we note that the boundary condition $\lambda=0$
at $\Sigma$, together with  \eqref{umbilic-e1}, implies
\begin{equation} \label{bdryeq1}
\nabla^2_g \lambda =-\frac1{n-1} g
\end{equation}
at $ \Sigma $.
Now choose a local orthonormal frames $ \{ e_i \} $ at the boundary
so that $e_i$ is tangential for $1 \le i \le n-1$ and $e_n = \nu$ is the
unit outward normal.  For $1\le i\le n-1$, \eqref{bdryeq1} implies
\begin{equation*}
 \begin{split}
  -\frac1{n-1}&= \nabla^2_g \lambda (e_i , e_i ) \\
  &= e_ie_i ( \lambda)  - \nabla_{e_i}e_i  (\lambda) \\
        & =-\langle \nabla_{e_i}e_i,e_n\rangle e_n(\lambda).
\end{split}
\end{equation*}
Summing over $1\le i\le n-1$, we have
\begin{equation}\label{umbilic-e2}
-1=H\frac{\p \lambda}{\p \nu},
\end{equation}
where  $H$ is the mean curvature of $ \Sigma $ with respect to $
\nu$. If $\lambda<0$ somewhere, then the set $U=\{x\in
\Omega|\lambda<0\}$ is a nonempty open set contained in $\Omega$
such that $\lambda=0$ at $ \p U$ because $\lambda=0$ at $\Sigma$.
Since $(n-1)\Delta \lambda +K\lambda =-n<0$, we have
$$
\int_U\lf((n-1)|\nabla \lambda|^2-K\lambda^2 \ri)dV < 0,
$$
contradicting the fact that the first Dirichlet eigenvalue of $(n-1)\Delta
+K$ is nonnegative. Therefore $\lambda>0$ in $\Omega$.
This together with \eqref{umbilic-e2} implies
$\frac{\p \lambda}{\p \nu}<0$ at $ \Sigma$.
Hence $H >0 $ and $e_n=\nu=-\nabla \lambda/|\nabla
\lambda|$.

Next, let $X$ and $ Y$ be two arbitrary vector fields on $ \Sigma$
which are tangential to $\Sigma$.   At $ \Sigma$, \eqref{bdryeq1}
implies
\begin{equation*}
\begin{split}
0  = & \nabla^2_g \lambda (X, \nu )  \\
    = & X \nu (\lambda)- \nabla_X \nu (\lambda)  \\
    = & - X ( | \nabla \lambda |)
\end{split}
\end{equation*}
where we used $\lambda=0$ at $\Sigma$,
and
\begin{equation*}
\begin{split}
- \frac1{n-1}g(X,Y) = & \nabla^2_g \lambda(X,Y) \\
= & X Y (\lambda) - \nabla_X Y ( \lambda ) \\
= &  \la \nabla_X Y , | \nabla \lambda | \nu \ra_g \\
= & - | \nabla \lambda | \Pi ,
\end{split}
\end{equation*}
where $ \Pi $ denotes the second fundamental form of $ \Sigma $
with respect to $ \nu$. Hence, along each connected component $
\Sigma_\alpha$ of $ \Sigma$, we conclude that  $ | \nabla \lambda
| $ is a positive constant and $ \Pi $ equals a positive constant
multiple of the induced metric on $ \Sigma $. (The constants may
depend on $ \alpha$.)

(iv) follows directly from the Gauss equation and (iii).
This completes the
proof of the theorem.
\end{proof}

Theorem \ref{critical-spaceforms-s2-t1} shows that
the standard metrics on  geodesic balls (which are
contained in a hemisphere  for the case of $\mathcal{S}^n$)
in space forms are critical points of the volume functional in the
corresponding spaces of metrics. It is natural to ask whether they are the
only critical points with  that  boundary condition.
Using Theorem \ref{umbilic-s2-t2}, we give some partial answer to
this question.

\begin{cor}\label{umbilic-s2-t2-c}
Let $ \Omega $ be an $ n $-dimensional
compact manifold with a smooth connected boundary
$ \Sigma $. Let $ \gamma $  be a given metric on $ \Sigma $
and let $ K = 0 $ or $ - n (n-1)$. Suppose
$ g \in \mathcal{M}^K_\gamma $ is a smooth metric
and $ g $ is a critical point of the volume functional
$ V (\cdot) $ in $  \mathcal{M}^K_\gamma $.
Let $ \nu $ be the outward unit normal vector to $ \Sigma$
in $ ( \Omega, g) $.
\begin{enumerate}

 \item [(i)] If $ K = 0 $, $ (\Sigma, \gamma) $ is
 isometric to a geodesic sphere in $ \R^n $
 and  $\Omega$ is  spin if $n\ge 8$,
 then $\Ric(\nu,\nu) $ is a non-positive constant along $ \Sigma $,
 and $\Ric(\nu,\nu) =0$ if and only if $(\Omega,g)$ is
 isometric to a standard  ball in $\R^n$.

    \item[(ii)]  If $ n = 3 $, $ K = 0 $, $ \Omega $ is oriented
    and $ Ric (\nu, \nu) = 0 $ along $ \Sigma $,
    then $(\Omega,g)$ is isometric to
    a standard ball in $\R^3$.

    \item[(iii)] If $ n = 3$, $ K  =  -6 $ and $ (\Sigma, \gamma) $ is
 isometric to a geodesic sphere in $ \H^3 $, then $ \Ric(\nu, \nu) $
    is a constant satisfying $\Ric(\nu,\nu)\le-2 $ along $ \Sigma$,
     and $ Ric ( \nu, \nu  ) = - 2 $ if and only
        if $(\Omega,g)$ is isometric to a geodesic ball in $\H^3$.
\end{enumerate}
\end{cor}

\begin{proof} (i)  As $ (\Sigma, \gamma) $ is isometric to a geodesic
sphere, say $ \Sigma_0 $,  in $ \R^n $, $ R_\Sigma $ is  a constant.
Hence, by (iii) and (iv) in Theorem \ref{umbilic-s2-t2},
$ Ric(\nu, \nu)$ is a constant along $ \Sigma $. On the other hand,
applying the Gauss equation to $ \Sigma_0 $ in $ \R^n$,
we have
\begin{equation} \label{s3-cor1-e2}
              R_\Sigma =  \frac{n-2}{n-1}H^2_0,
\end{equation}
where $ H_0 $ is the mean curvature of $ \Sigma_0 $ in $ \R^n $
(which is a constant). Hence, it follows from (iv) in Theorem \ref{umbilic-s2-t2},
\eqref{s3-cor1-e2} and the fact $ K = 0 $ that
\begin{equation}
2 \Ric(\nu,\nu)=\frac{n-2}{n-1}\lf[H^2 - H_0^2\ri].
\end{equation}
By the results in \cite{ShiTam02,Miao02}, which are
generalizations of the positive mass theorem in
\cite{SchoenYau79,Witten81,Schoen89}, we have $H\le H_0$, and $H=
H_0$ if and only if $(\Omega,g)$ is isometric to a standard ball
in $ \R^n $. From these, (i) follows.

(ii) By the assumption, $\Sigma $ is a connected orientable $
2$-surface. Let $ K_\Sigma $ be the Gaussian curvature of $
(\Sigma, \gamma) $. It follows from (iii), (iv) in Theorem
\ref{umbilic-s2-t2} and the facts $K=0$, $ Ric( \nu, \nu) = 0 $
that
$$ K_\Sigma = \frac{1}{4} H^2, $$
which is  a positive constant. Therefore,
$(\Sigma, \gamma )$ is isometric to a round sphere in
$ \R^3 $. Now (ii) follows from (i).

(iii) The proof is similar to the proof of (i).
 As $ (\Sigma, \gamma) $ is isometric to a geodesic
sphere, say $ \Sigma_1 $,  in $ \H^3 $, $ R_\Sigma $ is  a constant.
Hence, $ Ric(\nu, \nu)$ is a constant. On the other hand,
applying the Gauss equation to $ \Sigma_1 $ in $ \H^3$,
we have
\begin{equation} \label{s3-cor1-e3}
            2 (-2) +  R_\Sigma = -6 + \frac{n-2}{n-1}H^2_1,
\end{equation}
where $ H_1 $ is the mean curvature of $ \Sigma_1 $ in $ \H^3 $
(which is a constant). Hence, it follows from (iv) in Theorem \ref{umbilic-s2-t2},
\eqref{s3-cor1-e3} and the fact $ K = -6 $ that
\begin{equation}
2 [ \Ric(\nu,\nu) - ( -2) ] = \frac{n-2}{n-1}\lf[H^2 - H_1^2\ri].
\end{equation}
By Theorem 3.8 in \cite{ShiTam07}, we
we have $H \le H_1$, and $ H = H_1$ if and only if $(\Omega,g)$ is
isometric to  a geodesic ball in $ H^3 $. From these, (iii) follows.
\end{proof}

 As another application of Theorem \ref{umbilic-s2-t2},
 assuming $ (\Sigma, \gamma)$ can be
isometrically embedded as a convex hypersurface in
$ \mathbb{R}^n $, we can compare the
volume of any critical point in $ \mathcal{M}^0_\gamma $
with the Euclidean volume enclosed by $(\Sigma, \gamma)
$ in $\mathbb{R}^n$.

\begin{thm} \label{BY-volume-s2-t3}
Let $(\Omega, g)$ be an $ n$-dimensional smooth compact  Riemannian
manifold with zero scalar curvature, with a smooth
connected boundary $\Sigma$, such that $g$ is a critical point of
the volume functional in $\mathcal{M}_\gamma^0$
where $\gamma=g|_\Sigma$.
Suppose $(\Sigma, \gamma)$ can be isometrically embedded in $\R^n$
as a compact strictly convex hypersurface $\Sigma_0$.
If the dimension $ n \geq 8 $, $ \Omega $ is also assumed to be spin.
Then
$$  V(g) \ge V_0 $$
where $V(g)$ is the volume of $(\Omega,g)$ and $V_0$ is the
Euclidean volume of the domain bounded by $\Sigma_0$
in $ \mathbb{R}^n $.
Moreover, $ V (g) = V_0 $ if and
only if $(\Omega, g)$ is isometric to a standard ball in $\R^n$.
\end{thm}

\begin{proof} As $ g $ is a smooth critical point of the volume functional in
$ \mathcal{M}^0_\gamma$, by Theorem \ref{critical-nec-suff-s1t1}
there is a smooth function $ \lambda $ on $ \bar{\Omega} $ such that
$ \lambda = 0 $ at $ \Sigma$ and
\begin{equation*}
    -(\Delta_g \lambda)g+\nabla^2_g \lambda -\lambda \Ric(g) =g
\end{equation*}
in $ \Omega$.
As $ g $ has zero scalar curvature, the condition in (iii) in Theorem
 \ref{umbilic-s2-t2} is satisfied, therefore the mean curvature $H$
 of $\Sigma$ in $ (\Omega, g )$
 with respect to the unit outward normal $ \nu $ is a positive constant.
 Moreover, by \eqref{umbilic-e2}, $ H $ and $ \frac{ \p \lambda}{\p \nu}$
 satisfies
 $$ -1 = H \frac{ \p \lambda}{\p \nu} $$
 at $ \Sigma$. Integrating on $ \Sigma$,
 we have
\begin{equation}\label{BY-e1}
  \begin{split}
  |\Sigma|&=-\oint_{\Sigma} H\frac{\p \lambda}{\p \nu}\\
  &=-H\int_\Omega \Delta_g \lambda dV_g\\
  &=\frac{n}{n-1} H V(g) ,
\end{split}
\end{equation}
where the last step follows from the fact that
$$ \Delta_g \lambda = - \frac{ n }{n -1} . $$
On the other hand, by the results in \cite{ShiTam02},
\begin{equation}\label{BY-e2}
\oint_{\Sigma_0}H_0 \ge \oint_\Sigma H=|\Sigma|H
\end{equation}
where $H_0$ is the mean curvature of $\Sigma_0$ with respect to
the outward normal in $ \mathbb{R}^n$. By a Minkowski inequality
\cite{BonnesenFenchel}:
\begin{equation}\label{BY-e3}
|\Sigma|^2\ge \frac{n}{n-1}V_0\oint_{\Sigma_0}H_0.
\end{equation}
It follows from \eqref{BY-e1}, \eqref{BY-e2} and \eqref{BY-e3}
that
\begin{equation}\label{BY-e4}
\begin{split}
 \frac{|\Sigma|}{V(g)}&=\frac{n}{n-1} H \\
 &\le \frac{n}{(n-1)|\Sigma|}\oint_{\Sigma_0}H_0\\
 &\le \frac{|\Sigma|}{V_0}.
\end{split}
\end{equation}
Hence $V ( g) \geq V_0 $. If $ V(g) = V_0 $,  then \eqref{BY-e2}
becomes an equality. By the results in \cite{ShiTam02} again (see also \cite{Miao02}),
we know that $(\Omega,g)$ is isometric to a domain in $\R^n$. Finally
by Theorem \ref{critical-spaceforms-s2-t1}, we conclude that
$(\Omega, g)$ is isometric to a standard ball in $\mathbb{R}^n$.

\end{proof}

\section{Second variational formula for the volume functional} \label{sec-secondvar}
 In this section, we will compute the second variation
 of the volume functional $ V(\cdot) $ at critical points in
 $ \mathcal{M}^K_\gamma $. First, we give
 a formula for the second derivative of the scalar
 curvature. We remark that our notation convention for the curvature tensor
gives 
$$ R_{ijkl} = \kappa ( g_{ik} g_{jl} - g_{il} g_{jk} ) $$
in the case that $ g $ has constant sectional curvature $ \kappa
$.

\begin{lma}\label{s3-l1}
Let $ \{ g(t) \} $ be a smooth path of $ C^2 $ metrics with $ g(0) = g$.
Let $ R(t) $ be the scalar curvature of $ g(t) $. Then
\begin{equation} \label{2variationR}
\begin{split}
R^{\prime \prime}(0) = & \ \Delta_g ( | h |_g^2 )
+ 2 \la h, \nabla^2_g ( \tr_g h ) \ra_g
- 4  \la \nabla_g  ( \div_g h) , h\ra_g   \\
&  - 2  | \div_g  h   - \frac{1}{2} d ( \tr_g h ) |_g^2
- \frac{1}{2} | \nabla_g h |_g^2    - g^{pq} h^{lk}_{ \ \ ;p} h_{lq;k} \\
&
 + 2 h^{lp} R_{lkps} h^{sk} + D\mathcal{R}_{g} ( h^\prime) ,
\end{split}
\end{equation}
where $ h = g^\prime ( 0 ) $, $ h^\prime = g^{ \prime \prime }(0)
$, $ \nabla_g (\cdot)$ and  ` $;$ ' denote covariant derivative
with respect to $ g $, $ R_{ijkl} $ is the curvature tensor of $ g
$, and $ D\mathcal{R}_{g}( \cdot ) $ is the linearization of the
scalar curvature map $ \mathcal{R} $ at $ g $.
\end{lma}

\begin{proof} The first derivative of the scalar curvature is given by
\begin{equation}\label{scalar-2nd-e2}
    R^{\prime} (t) =-\Delta_g  ( \tr_g g^\prime ) +\div_g (\div_g  ( g^\prime) )
    -\la g^\prime , \Ric \ra_g .
\end{equation}
For any fixed point $ p $, let $ \{ x_i \}$ be a normal coordinate
chart at $ p $ with respect to $ g(0) = g $. We will use
  `$  , $'   to denote partial derivative and  `$ ;$'  to
  denote covariant derivative. At $ t = 0 $, we have
\begin{equation}\label{scalar-2nd-e3}
\begin{split}
    - [ \Delta_g ( \tr_g g^\prime ) ]'  = & -
    \lf[ g^{ij} \lf(\tr_g g^\prime )_{,ij} - g^{ij} \Gamma_{ij}^k
    ( \tr_g g^\prime)_{,k}\ri)\ri]'  \\
       =  &-\bigg[-g^{ip}g^{qj}h_{pq}\lf(\tr_g h \ri)_{,ij}-\Delta_g |h|_g^2+
        \Delta_g ( \tr_g h' )\\
       & \quad-\frac12g^{ij}g^{kl}
        \lf(h_{jl;i}+h_{il;j}-h_{ij;l}\ri)\lf(\tr_g h \ri)_{,k} \bigg]\\
       = & \ \la h, \nabla^2_g ( \tr_g h ) \ra_g +\Delta_g |h|_g^2-
        \Delta_g ( \tr_g h' ) \\
        &  + \la \div_g\, h, d ( \tr_g h ) \ra_g - \frac12 |d ( \tr_g h ) |_g^2,
\end{split}
\end{equation}
\begin{equation}\label{scalar-2nd-e4}
\begin{split}
\lf[ \div_g \lf( \div_g \, g^\prime \ri) \ri]'  &=
\lf( g^{ij}g^{kl}\ri)'h_{jl;ki}+g^{ij} g^{kl}\lf(g^\prime_{jl;ki}  \ri)' \\
&=-\la \nabla_g\div_gh,h\ra_g-g^{ij}g^{kp}g^{ql}h_{pq}h_{jl;ki}
+g^{ij} g^{kl}\lf( g^\prime_{jl;ki} \ri)' \\
&=-2\la \nabla_g\div_gh,h\ra_g+h^{lp}R_{lkps}h^{sk}-h^{ij}g^{kl}R_{jk}h_{il} \\
& \ \ \ \ +g^{ij} g^{kl}\lf( g^\prime_{jl;ki} \ri)' ,
\end{split}
\end{equation}
where we have used the Ricci identities. Now
\begin{equation*}\label{scalar-2nd-e5}
\begin{split}
 g^{ij} g^{kl}\lf( g^\prime_{jl;ki} \ri)'
= &  g^{ij}g^{kl}\lf[ (g^\prime_{jl;k} )_{,i}
-\Gamma_{ij}^s g^\prime_{sl;k} -\Gamma_{il}^s g^\prime_{js;k}
-\Gamma_{ik}^s g^\prime_{jl;s}\ri]' \\
= & g^{ij}g^{kl}\bigg[
\lf( h'_{jl;k} -\lf(\Gamma_{kj}^s\ri)'h_{sl}-\lf(\Gamma_{kl}^s\ri)'h_{js}\ri)_{,i}
\\
&  -\lf(\Gamma_{ij}^s\ri)'h_{sl;k}
-\lf(\Gamma_{il}^s\ri)'h_{js;k}
 -\lf(\Gamma_{ik}^s\ri)'h_{jl;s}\bigg],
\end{split}
\end{equation*}
and
\begin{equation*}
g^{ij}g^{kl} \lf(h'_{jl;k}\ri)_{,i} = \div_g (\div_g  h') ,
\end{equation*}
\begin{equation*}
\begin{split}
-g^{ij}g^{kl}
\lf(\lf(\Gamma_{kj}^s\ri)'h_{sl}\ri)_{,i} &=-\frac12 g^{ij}g^{kl}
\lf(g^{sm}\lf(h_{km;j}+h_{mj;k}-h_{kj;m}\ri)h_{sl}\ri)_{,i}\\
&=-\frac12 g^{ij}g^{kl}
\lf(g^{sm} h_{km;j}h_{sl}\ri)_{,i}\\
&=-\frac14 \Delta_g |h|_g^2 ,
\end{split}
\end{equation*}
\begin{equation*}
\begin{split}
-g^{ij}g^{kl} \lf(\lf(\Gamma_{kl}^s\ri)'h_{js}\ri)_{,i}
= & -\frac12 g^{ij}g^{kl}
\lf(g^{sm}\lf(h_{km;l}+h_{ml;k}-h_{kl;m}\ri)h_{js}\ri)_{,i}\\
= & -\la\nabla_g \div_g  h, h\ra_g
+ \frac12 \la \nabla^2_g ( \tr_g h ), h \ra_g \\
&  -|\div_g h|_g^2 +\frac12 \la \div_g h, d ( \tr_g h ) \ra_g ,
\end{split}
\end{equation*}
\begin{equation*}
\begin{split}
-g^{ij}g^{kl}\lf(\Gamma_{ij}^s\ri)'h_{sl;k}&=-\frac12g^{ij}g^{kl}g^{sm}
\lf(h_{im;j}+h_{mj;i}-h_{ij;m}\ri)h_{sl;k}\\
&= - |\div_g h|_g^2 + \frac12 \la \div_g h , d ( \tr_g h ) \ra_g ,
\end{split}
\end{equation*}
\begin{equation*}
\begin{split}
-g^{ij}g^{kl}\lf(\Gamma_{il}^s\ri)'h_{js;k}&=-\frac12g^{ij}g^{kl}g^{sm}
\lf(h_{im;l}+h_{ml;i}-h_{il;m}\ri)h_{js;k}\\
&=-\frac12 |\nabla_g h|_g^2 ,
\end{split}
\end{equation*}
\begin{equation*}
\begin{split}
-g^{ij}g^{kl}\lf(\Gamma_{ik}^s\ri)'h_{jl;s}&=-\frac12g^{ij}g^{kl}g^{sm}
\lf(h_{im;k}+h_{mk;i}-h_{ik;m}\ri)h_{jl;s}\\
&=- g^{pq}h^{lk}_{\ \ ;p}h_{lq;k} + \frac12 |\nabla_g h|_g^2.
 \end{split}
\end{equation*}
Hence,
\begin{equation}\label{scalar-2nd-e6}
\begin{split}
\lf( \div_g \lf(\div_g  g^\prime \ri)\ri)'
= & -3\la \nabla_g\div_gh,h\ra_g+h^{lp}R_{lkps}h^{sk}
-h^{ij}g^{kl}R_{jk}h_{il} \\
& +\div_g(\div_g\, h') - \frac14\Delta_g |h|_g^2
  +\frac12 \la \nabla^2_g ( \tr_g h ), h\ra_g \\
 &  + \la\div_g h,d ( \tr_g h) \ra_g-2|\div_g\, h|_g^2- g^{pq}h^{lk}_{\ \ ;p}h_{lq;k}.
\end{split}
\end{equation}
Next,
\begin{equation}
\begin{split}
- [ \la g^\prime , \Ric \ra_g ]^\prime
&=-\lf(g^{ij}g^{kl} g^\prime_{ik}R_{jl}\ri)'\\
&= 2h^{ij}g^{kl}R_{jk}h_{il} - \la h',\Ric\ra_g -
g^{ij}g^{kl}h_{ik}R_{jl}^\prime ,
\end{split}
\end{equation}
\begin{equation*}
\begin{split}
g^{ij}g^{kl}h_{ik}R_{jl}' = & \
\frac12g^{ij}g^{kl}h_{ik}g^{pq}\lf(h_{lp;jq}+h_{jp;lq}-h_{jl;pq}-h_{pq;jl}\ri) \\
= & -\frac14 \Delta_g ( |h|_g^2) + \frac12 |\nabla_g h|_g^2-
\frac12 \la \nabla^2_g (\tr_g h) , h \ra_g \\
&  +g^{ij}g^{kl}g^{pq}h_{ik}h_{lp;jq}\\
 =&-\frac14 \Delta_g ( |h|_g^2) + \frac12 |\nabla_g h|_g^2-
\frac12 \la \nabla^2_g (\tr_g h) , h \ra_g \\
&+\la\nabla_g\div_gh,h\ra_g+h^{ij}g^{kl}R_{jk}h_{il}-h^{lp} R_{lkps} h^{sk},
\end{split}
\end{equation*}
where we have used the Ricci identity in the last step. Hence
\begin{equation}\label{scalar-2nd-e7}
\begin{split}
 - [ \la g^\prime , \Ric \ra_g ]^\prime
& =   h^{ij}g^{kl}R_{jk}h_{il}-\la h',\Ric\ra_g
+\frac14 \Delta_g |h|_g^2 -\frac12|\nabla_g h|^2\\
&+\frac12 \la \nabla^2_g ( \tr_g h ), h\ra_g -\la\nabla_g\div_gh,h\ra_g
+h^{lp} R_{lkps} h^{sk}.
\end{split}
\end{equation}
So  \eqref{2variationR} follows from \eqref{scalar-2nd-e2},
\eqref{scalar-2nd-e3}, \eqref{scalar-2nd-e6} and
\eqref{scalar-2nd-e7}.
\end{proof}

Now we are in a position to compute the second variation of 
$ V(\cdot)$ in $ \mathcal{M}^K_\gamma$. We state the 
formula in a general setting which does not require the 
manifold structure of $ \mathcal{M}^K_\gamma$.

\begin{thm} \label{svar-formula}
Let $ \Omega $ be an $ n$-dimensional connected compact 
manifold with a smooth boundary $ \Sigma $. Suppose $ g $ 
is a smooth metric on $ \Omega $ such that there is a 
smooth function $ \lambda $ on $ \Omega $ satisfying
\begin{equation} \label{sec4-critical-e1}
   \left\{
\begin{array}{ll}
    -(\Delta_g\lambda)g+\nabla^2_g\lambda-\lambda \Ric(g)
    &=g \hbox{\ in $\Omega$ ;} \\
    \lambda &=0 \hbox{\ on $\Sigma$.}
\end{array}
\right.
\end{equation}
Let $ \gamma = g |_{T(\Sigma)} $ and let $ K $ be the constant
that equals the scalar curvature of $ g $. Suppose $ \{ g(t) \} $
is a smooth path of metrics
in $ \mathcal{M}_\gamma $  with $ g(0) = g $ and
$ g(t) \in \mathcal{M}_\gamma^K $.
Let $ V(t) $ be the volume of $ ( \Omega, g(t) ) $, then
\begin{equation} \label{svar}
\begin{split}
V^{\prime \prime} (0)  = &
 \int_\Omega  \lf\{  \frac{1}{4}  ( \tr_g h )^2
 +   \l  \lf[  | \div_g ( h )  - \frac{1}{2} d ( \tr_g h ) |_g^2
+ \frac{1}{4}  | \nabla_g h |_g^2  \ri]  \ri\} dV_g \\
& + \int_\Omega   \l \lf[ \la \nabla_g \div_g h, h\ra_g
-   \la h ,  \nabla^2_g ( \tr_g h )  \ra_g \ri] d V_g \\
 &  - \int_\Omega   \frac{1}{2} \l
 \lf[ | \div_g h |_g^2 + h^{sp} R_{kpls} h^{lk} \ri] d V_g ,
\end{split}
\end{equation}
where $ h = g^\prime(0) $.
\end{thm}

\begin{proof}
We first note that $ V^\prime ( 0 ) = 0 $ by (\ref{sec4-critical-e1})
and the proof of Theorem \ref{critical-nec-suff-s1t1}.
To compute $ V^{\prime \prime }(0)$, we start with
\begin{equation} \label{vdoublep}
V^{\prime \prime}(0)  =  \int_\Omega
\lf[ \frac14 \left( \tr_g h \right)^2 - \frac12 |h |_g^2 \ +
 \frac12  ( \tr_g h^{ \prime} ) \ri] \ d V_g ,
\end{equation}
where $ h = g^\prime(0) $ and $ h^\prime = g^{\prime \prime}(0) $.
Our aim is to express the last integral in terms of $h$.

Recall that (\ref{sec4-critical-e1}) implies \be \label{lap-lambda} \Delta_g
\l = - \frac{1}{n-1} ( K   \l + n) \ee and \be \label{hes-lambda}
\nabla_g^2 \l = \l \Ric - \frac{1}{n-1} \lf( \l K + 1 \ri) g  .
\ee
Applying  the fact that $ g(t) $ has constant scalar curvature $ K
$ and Lemma \ref{s3-l1}, we have
\begin{equation} \label{svar-C}
\begin{split}
0 = & R^{\prime \prime} (0) \\
 =  & I -  \Delta_g ( \tr_g h^\prime ) + \div_g ( \div_g h^\prime) -
 \la h^\prime, \Ric \ra_g  ,
\end{split}
\end{equation}
where
\begin{equation} \label{defofI}
\begin{split}
I = & \Delta_g ( | h |_g^2 )
+ 2 \la h, \nabla^2_g ( \tr_g h ) \ra_g
- 4  \la \nabla_g  \div_g h , h\ra_g   \\
&  - 2  | \div_g ( h )  - \frac{1}{2} d ( \tr_g h ) |_g^2
- \frac{1}{2} | \nabla_g h |_g^2    - g^{pq} h^{lk}_{ \ \ ;p} h_{lq;k} \\
&
 + 2 h^{lp} R_{lkps} h^{sk}.
 \end{split}
\end{equation}
Note that $I$ involves only $h$ and its derivatives. In what
follows, we will omit the volume form and the area form in
integrals  for convenience. All integrals are taken with
respect to the metric $g = g(0) $.
Integrating by parts, we have
\begin{equation} \label{svar-A}
\begin{split}
& \ \int_\Omega  \lf[-\l \Delta_g(\tr_g(h'))+ \l \div_g ( \div_g h^{ \prime} )
-  \l \la h^\prime , \Ric \ra_g \ri]  \\
 = & \
\int_\Omega \lf[ -( \tr_g h^{ \prime} ) ( \Delta_g \l ) + \la \nabla_g^2 \l , h^{ \prime} \ra_g -  \l \la h^\prime , \Ric \ra_g \ri] \\
& \ + \oint_{\Sigma} \lf[ \lambda_{;\nu} \lf( \tr_g h^{ \prime} \ri)
-  h^\prime (\nu, \nabla \l) \ri]
\\
 = & \ \int_\Omega     \tr_g  h^{ \prime} +\oint_{\Sigma} \lambda_{;\nu} \lf( \tr_g h^{ \prime} \ri)
- \oint_{\Sigma}  h^\prime (\nu, \nabla \l)
\end{split}
\end{equation}
where we have used $\lambda |_\Sigma =0$,
\eqref{lap-lambda} and \eqref{hes-lambda}.
The condition $\lambda |_\Sigma =0 $, together with
$ h^\prime |_{ T(\Sigma) } = 0 $, also implies
\begin{equation*}
   -  \lambda_{;\nu} ( \tr_g h^{ \prime} )
   + h^{\prime} (\nu, \nabla \l ) = 0  \ \mathrm{on} \ \Sigma.
\end{equation*}
Combining this with  \eqref{svar-C} and \eqref{svar-A}, we have,
\be \label{simpletrace}
 \int_\Omega ( \tr_g h^{\prime} )  = - \int_\Omega \l   I  .
\ee

Next we compute the integral of $\l I$. Integrating by
parts and applying (\ref{lap-lambda}) and (\ref{hes-lambda}), we
have
\begin{equation} \label{svar-J}
\begin{split}
\int_\Omega \l \Delta_g ( | h |_g^2 )
= &  \int_\Omega ( \Delta_g \l ) | h|_g^2
- \oint_{\Sigma} \l_{;\nu}  | h|_g^2 \\
= & -  \int_\Omega \frac{ 1 }{ n-1 } ( K \l + n)  | h|_g^2
- \oint_{\Sigma} \l_{;\nu}  |h|_g^2 .
\end{split}
\end{equation}
Integrating by parts, we also have
\begin{equation} \label{svar-F}
\begin{split}
\int_\Omega \l  g^{pq} h^{lk}_{ \ \ ;p} h_{lq;k}  = &
\int_\Omega ( \l g^{pq} h^{lk} h_{lq;k} )_{;p} -
g^{pq} \l_{;p} h^{lk} h_{lq;k}
-  \l g^{pq} h^{lk} h_{lq;kp}  \\
 = &   \int_\Omega - g^{pq} \l_{;p} h^{lk} h_{lq;k} -
 \l g^{pq} h^{lk} h_{lq;kp} .
\end{split}
\end{equation}
Now
\begin{equation} \label{svar-G}
\begin{split}
 \int_\Omega & g^{pq} \l_{;p} h^{lk} h_{lq;k} \\
   = &
 \int_\Omega ( g^{pq} \l_{;p} h^{lk} h_{lq} )_{;k} - g^{pq}
 \l_{;pk} h^{lk} h_{lq}   - g^{pq} \l_{;p} h^{lk}_{\ \ ;k}  h_{lq} \\
 = & \oint_{\Sigma} g^{pq} \l_{;p} h^{l \nu} h_{lq} +  \int_\Omega
  \frac{1}{ n -1 }  \left( K \l  + 1 \right) | h |_g^2
   -  \int_\Omega  \l  g^{pq} R_{pk} h^{lk}h_{lq} - \int_\Omega \l_{;p}
 (\div_g h)_l  h^{lp}\\
 =& \oint_{\Sigma}   \l_{;\nu}|h|^2_g  +  \int_\Omega
  \frac{1}{ n -1 }  \left( K \l  + 1 \right) | h |_g^2   -  \int_\Omega  \l  g^{pq} R_{pk} h^{lk}h_{lq} - \int_\Omega \l_{;p}
 (\div_g h)_l  h^{lp}
\end{split}
\end{equation}
where we have used the fact that $h|_{T(\Sigma)}=0$. Also
\begin{equation} \label{svar-H}
\begin{split}
\int_\Omega  \l_{;p} (\div_g h)_l  h^{lp}  = &
\int_\Omega [ \l ( \div_g h)_l  h^{lp} ]_{;p}   -  \l ( \div_g h)_{l;p} h^{lp}
-  \l | \div_g h |^2_g \\
 = &  - \int_\Omega   \l \la \nabla_g \div_g h, h\ra_g
 - \int_\Omega \l | \div_g h |_g^2 .
\end{split}
\end{equation}
By the Ricci identify,
\begin{equation} \label{svar-I}
\begin{split}
\int_\Omega
 \l g^{pq} h^{lk} h_{lq;kp} = & \int_\Omega  \l \la \nabla_g \div_g h, h\ra_g+
 \l g^{pq}R_{pk}h^{lk}h_{lq}-\l h^{sp}R_{kpls}h^{lk} .
 \end{split}
\end{equation}
 It follows from (\ref{svar-F})--\eqref{svar-I}  that
\begin{equation} \label{svar-K}
\begin{split}
\int_\Omega \l  g^{pq} h^{lk}_{ \ \ ;p} h_{lq;k}
= & - \oint_{\Sigma}    \l_{;\nu} |h|^2_g -  \int_\Omega
  \frac{1}{ n -1 }  \left( K \l  + 1 \right) | h |_g^2 \\
 &  - \int_\Omega 2 \l \la\nabla \div_g h,h\ra_g - \int_\Omega \l | \div_g h |_g^2 \\
 & + \int_\Omega \l h^{sp} R_{kpls} h^{lk}  .
\end{split}
\end{equation}
Combining \eqref{defofI}, (\ref{svar-J}) and (\ref{svar-K})   we
have
\begin{equation} \label{svar-L}
\begin{split}
  \int_\Omega \l I
= &   \    \int_\Omega - | h|_g^2 +2\l \la h,\nabla^2_g(\tr_g(h))\ra_g -2\l \la\nabla\div_g(h),h\ra_g\\
& - \int_\Omega  \lf[ 2 \l  | \div_g ( h )  - \frac{1}{2} d (
\tr_g h ) |_g^2
+ \frac{1}{2} \l | \nabla h |_g^2  \ri] \\
 &  + \int_\Omega  \l | \div_g h|_g^2 + \l h^{sp} R_{kpls} h^{lk}.
\end{split}
\end{equation}
Now (\ref{svar}) follows from (\ref{vdoublep}), (\ref{simpletrace})
and (\ref{svar-L}).
\end{proof}

Next we focus on domains $\Omega $ in space forms on
which the standard metrics are critical points of $ V(\cdot) $.
By Theorem \ref{critical-spaceforms-s2-t1},
those domains are precisely geodesic balls.

\begin{thm} \label{svar-positive}
Let $ \Omega $ be a geodesic ball with compact closure
in $ \mathbb{R}^n $, $\mathbb{H}^n $ or $ \mathbb{S}^n_+$.
Let $ \Sigma $ be  its boundary and $ g $ be the standard metric
on $ \Omega $. Let  $ \gamma = g |_{ T (\Sigma) }$
 and $ K $ be the constant that equals the scalar curvature of $ g $.
Let $ \{ g(t) \}$ be a smooth path of metrics in
$ \mathcal{M}^K_\gamma $ with  $ g(0) = g $. Let $ h $
denote $ g^\prime(0) $ which is assumed to be nonzero.

\begin{enumerate}
\item[(i)] If $ \Omega \subset \mathbb{R}^n $ or
$ \Omega \subset  \mathbb{S}^n_+$, then
$$ \frac{d^2}{dt^2} V ( g(t) ) |_{t=0} > 0 $$
for any  $ h $ satisfying
$$ \div_g h = 0 \ \ \mathrm{and} \ \ \tr_g h = 0 .$$

\item[(ii)] If $ \Omega \subset \mathbb{H}^n $, then
for any $ p \in \Omega $, there is a constant $\delta > 0 $,
depending on $ p $ and $ (\Omega, g)$, such that
$$ \frac{d^2}{dt^2} V ( g(t) ) |_{t=0} > 0 $$
for  any $ h $ which has compact support  in
$ B_\delta (p) \subset \Omega$ and satisfies
$$ \div_g h = 0 \ \ \mathrm{and} \ \ \tr_g h = 0 . $$
Here $ B_\delta (p) $ denotes the geodesic ball with
radius $ \delta $ centered at $ p $ in $ \mathbb{H}^n $.
\end{enumerate}
\end{thm}

\begin{proof}
From the proof of Theorem \ref{critical-spaceforms-s2-t1},
we know there exists a smooth positive function $ \lambda $
on $ \bar{\Omega} $ satisfying \eqref{sec4-critical-e1}.
(The explicit expression of $ \l $ is given by \eqref{linR}, 
\eqref{linH} and \eqref{linS}).  Applying Theorem
\ref{svar-formula} and the assumption $ \div_g h = 0 $ and $ \tr_g
h = 0$, we have
\begin{equation}  \label{svar-spaceform}
\begin{split}
\frac{d^2}{dt^2} V(t) |_{t = 0}
 =  \frac{1}{4}  \int_\Omega  \l    | \nabla_g h |_g^2
  -   \frac{1}{2} \int_\Omega  \l   h^{sp} R_{kpls} h^{lk} .
\end{split}
\end{equation}
As $ g $ has constant sectional curvature $ \frac{K}{n ( n-1)}$,
the integrant in the last integral above reduces to
\begin{equation*}
\begin{split}
   \lambda h^{sp} R_{kpls} h^{lk}  =  & \frac{K \lambda }{n ( n -1 )}
\lf[ ( \tr_g h )^2 - | h |_g^2 \ri]  .
\end{split}
\end{equation*}
Therefore (\ref{svar-spaceform}) becomes
\begin{equation} \label{svar-spaceform-reduced}
\begin{split}
\frac{d^2}{dt^2} V (t) |_{t = 0}
 =  \frac{1}{4}  \int_\Omega  \l    | \nabla_g h |_g^2
  +  \frac{ K }{2 n ( n -1 ) } \int_\Omega  \l  | h |^2_g  ,
\end{split}
\end{equation}
where we again used the assumption $ \tr_g h = 0 $.
As $ \l > 0 $ on $ \Omega$, (i) follows directly from
(\ref{svar-spaceform-reduced}).

Now suppose $ \Omega \subset \mathbb{H}^n_+ $,
we then have $ K = - n ( n -1) $ and
\begin{equation} \label{svar-H-reduced}
\begin{split}
\frac{d^2}{dt^2} V (t) |_{t = 0}
 =  \frac{1}{4}  \int_\Omega  \l    | \nabla_g h |_g^2
 -  \frac{ 1 }{2  } \int_\Omega  \l  | h |^2_g  .
\end{split}
\end{equation}
For any $ p \in \Omega $, first choose $ \delta > 0 $ such that
$ B_{\delta } (p) \subset \Omega $ and
$$ \min_{B_{\delta } (p)} \l \geq \frac{1}{2} \l (p)  , \ \
\max_{B_{\delta } (p)} \l \leq 2 \l (p) . $$
If $ h $ has compact support contained in
$ B_\delta  = B_{\delta} (p) $, then
\begin{equation}
\begin{split}
\frac{d^2}{dt^2} V (t) |_{t = 0}
 \geq & \l (p)  \lf\{ \frac{ 1}{8}  \int_{B_\delta}  | \nabla_g h |_g^2
 -   \int_{B_\delta}  | h |^2_g  \ri\}  \\
  \geq  & \l (p)  \lf\{ \frac{ 1}{8}  \int_{ B_\delta }
 | \nabla | h |_g  |_g^2
 -   \int_{B_\delta}  | h |^2_g  \ri\},
\end{split}
\end{equation}
where we used  $ | \nabla_g h |_g \geq | \nabla | h|_g |_g $
in the second step.  Let $ \l_1 (B) $ be the first eigenvalue
of $ \Delta_g $ on $ B_\delta $, then
\begin{equation}
\begin{split}
\frac{d^2}{dt^2} V (t) |_{t = 0}
  \geq  & \l (p)  \lf\{ \frac{ 1}{8}  \l_1 ( B_\delta ) - 1 \ri\}
   \int_{ B_\delta}    | h |^2_g   .
\end{split}
\end{equation}
 As
$ \lim_{\delta \rightarrow 0} \l_1 (B_\delta) = + \infty , $
we conclude that $ \frac{d^2}{dt^2} V (t) |_{t = 0} > 0
$ if $ \delta $ is smaller than some constant $\delta_0 $
depending only on $ p $ and $ (\Omega, g)$.
Therefore, (ii) is proved.
\end{proof}

\begin{lma} \label{constructg(t)}
Suppose $g\in \mathcal{M}_\gamma^K $ is a smooth Einstein
metric satisfying the property that the first eigenvalue of
$(n-1)\Delta_g+K$ is positive. Let $h$ be a symmetric
(0,2) tensor on $\Omega$ such that $h |_{ T(\Sigma)} = 0 $,
$ \tr_g h = 0 $ and $ \div_g \div_g(h)=0$.
Then there is a variation $ \{ g(t) \}  \subset
\mathcal{M}_\gamma^K$ of $ g $ such that  $g'(0)=h$.
\end{lma}
\begin{proof} As in the proof of Theorem \ref{critical-nec-suff-s1t1},
  we can find $g(t) \in \mathcal{M}_\gamma^K$ such that $g(t)$
  is of the form $u(t)^{\frac{4}{n-2}} (g+th)$ with $u(0)=1$. Moreover
  $ v = u'(0)$ satisfies:
\begin{equation}
  \begin{cases}
    (n-1)\Delta_{g} v +Kv&= \lf( \frac{n-2}{4} \ri)
    D \mathcal{R}_g ( h ), \text{\ in $\Omega$}\\
    v&=0,  \ \ \ \ \ \ \ \ \ \  \ \ \ \ \ \ \ \text{\ on $\Sigma$.}
\end{cases}
    \end{equation}
    As $\tr_g(h)=0$,  $\div_g(\div_g(h))=0$ and $g$ is Einstein,
    we know $ D \mathcal{R}_g ( h ) = 0 $
    by \eqref{scalar-2nd-e2}. As the first eigenvalue of
   $(n-1)\Delta_g+K$ is positive, we have $v  \equiv0$, hence
   $g'(0)=h$.
\end{proof}

In the appendix following this section, we will explicitly construct
 trace free and divergence free  $(0,2)$ symmetric tensors
with prescribed compact support on space forms (see also
\cite{Beig97,Corvino07}). Therefore, by the existence of such
tensors and by Theorem \ref{svar-positive} and Lemma
\ref{constructg(t)}, we have:

\begin{cor} Let $(\Omega,g)$ be given as in
Theorem \ref{svar-positive}. There exists a
 variation $ \{ g(t) \}  \subset \mathcal{M}_\gamma^K$
of $ g $  such that $ \frac{d}{dt} V ( g(t) ) |_{t = 0} = 0 $ and
$$ \frac{d^2}{dt^2} V ( g(t) ) |_{t=0} > 0 . $$
In particular, the volume of the standard metric $ g $
is a strict local minimum along such a variation. 
\end{cor}

Next, we show that there indeed exist deformations along 
which the volume of the standard metric is a strict local maximum.

\begin{thm} \label{svar-negative}
Let $ \Omega $ be a geodesic ball with compact closure
in $ \mathbb{R}^n $, $\mathbb{H}^n $ or $ \mathbb{S}^n_+$.
Let $ \Sigma $ be  its boundary and $ g $ be the standard metric.
Let  $ \gamma = g |_{ T (\Sigma) }$
 and $ K $ be the constant that equals the scalar curvature of $ g $.
Suppose the dimension $ n $ satisfies $ 3 \leq n \leq 5 $.
\begin{enumerate}

\item[(i)] If $ \Omega \subset \mathbb{R}^n $,
there exists a variation $ \{ g(t) \}  \subset \mathcal{M}_\gamma^K$
of $ g $  such that $ \frac{d}{dt} V ( g(t) ) |_{t = 0} = 0 $ and
$$ \frac{d^2}{dt^2} V ( g(t) ) |_{t=0}  < 0 . $$

\item[(ii)] If $ \Omega \subset \mathbb{H}^n $ or $ \mathbb{S}^n_+$,
there exists a small positive constant $ \delta $
such that, if  the geodesic radius of $ \Omega $ is less than $ \delta $,
then there exists a variation $ \{ g(t) \}  \subset \mathcal{M}_\gamma^K$
of $ g $  such that $ \frac{d}{dt} V ( g(t) ) |_{t = 0} = 0 $ and
$$ \frac{d^2}{dt^2} V ( g(t) ) |_{t=0}  < 0 . $$

\end{enumerate}
\end{thm}

\begin{proof}
We first consider the case $ \Omega \subset \mathbb{R}^n$. We use
$ g_0 $ to denote the standard Euclidean metric on $\mathbb{R}^n$
and $ \nabla_0 ( \cdot) $ to denote the covariant derivative taken
with respect to $ g_0$. Let $ \hh $ be an arbitrary, nonzero, symmetric
$(0,2)$ tensor that is parallel on  $(\Omega, g_0)$ and satisfies
$\tr_{g_0} \hh = 0 $. Define $ h = \l \hh .$ Then  
$ h |_{T(\Sigma) } = 0 $ and satisfies
\begin{equation} \label{svar-M}
\tr_{g_0} h = 0 \ \
\end{equation}
\begin{equation} \label{svar-N}
(\div_{g_0} h)_{i}   =  g^{jk} \l_{;k} \hh_{ij} ,
\end{equation}
\begin{equation} \label{svar-O}
\div_{g_0} ( \div_{g_0}  h) =
\la \nabla^2_{g_0} \l , \hh \ra_{g_0} = 0 ,
\end{equation}
where we have used the equation
\begin{equation} \label{svar-P}
\nabla^2_{g_0} \l = - \frac{1}{n-1} g_0 .
\end{equation}
By Lemma \ref{constructg(t)}, $ h \in T_{g_0} \mathcal{M}^K_\gamma
$.   Plug this $ h $ into (\ref{svar}) in Theorem
\ref{svar-formula} and use the fact $g_0$ has zero curvature, we
have
\begin{equation*} \label{svar-Q}
\begin{split}
\frac{d^2}{dt^2} V ( g(t) ) |_{t =  0}
 = &   \int_\Omega  \l  \lf[ \frac{1}{2} | \div_{g_0} ( h )   |_g^2
+ \frac{1}{4}  | \nabla_0 h |_{g_0}^2
+ \la \nabla_0 ( \div_{g_0} h ), h \ra_{g_0}  \ri] .
\end{split}
\end{equation*}
Applying (\ref{svar-N}) and (\ref{svar-P}), integrating by parts
and using the fact that $ \hh $ is parallel, we have
\begin{equation}
\begin{split}
\int_\Omega \l | \div_{g_0}  h |^2_{g_0}
= & \frac{1}{2(n-1)} \int_\Omega \l^2 | \hh |_{g_0}^2 .
\end{split}
\end{equation}
\begin{equation}
\begin{split}
\int_\Omega \l | \nabla h |^2_{g_0}
= & \frac{n}{2(n-1)}  \int_\Omega \l^2  | \hh |_{g_0}^2 ,
\end{split}
\end{equation}
\begin{equation}
\begin{split}
\int_\Omega \l ( \div_{g_0} h )_{l;p} h^{lp}
= & - \frac{1}{n-1} \int_\Omega  \l^2 | \hh |^2_{g_0} \\
\end{split}
\end{equation}
Thus,
\begin{equation} \label{svar-R}
\begin{split}
V^{\prime \prime} (0)
 =  & \frac18 \lf( \frac{ n - 6 }{  n - 1 } \ri) \int_\Omega \l^2  | \hh |^2_{g_0}
 < 0 ,
 \end{split}
\end{equation}
for $ n = 3, 4, 5 $.

Next, we consider the case $ \Omega \subset \mathbb{H}^n$.
For any $ \kappa > 0 $, consider the metric
$ g_\kappa =  \lf( 1 - \frac{ \kappa^2 }{ 4 } | x |^2 \ri)^{ - 2 } g_0 , $
which is defined on $ \{ | x | < \frac{ 2 }{  \kappa } \} $
 and has constant sectional curvature $ - \kappa^2 $.
Let $ B = \{ | x | < 1 \} $ and $ \Sigma $ be its boundary.
Let $ \l_k $ be the smooth
function on $ B $ defined by
\begin{equation}
\l_k = \frac{ 1 }{ ( n-1 ) \kappa^2 }
\lf( 1 - \frac{ 4 - \kappa^2 }{ 4 + \kappa^2 } \cosh{ \kappa r } \ri),
\end{equation}
where $ r $ is the geodesic distance from the origin in
$( B, g_\kappa )$. We also define
\begin{equation}
\l_0 = \frac{ 1 }{ 2 ( n -1)  } \lf( 1 - | x |^2 \ri).
\end{equation}
For $ \kappa \geq 0 $, the function  $ \l_\kappa $ satisifes
\begin{equation}
\lf\{
\begin{array}{lcc}
- ( \Delta_{ g_\kappa } \l_\kappa ) g_\kappa
+ \nabla^2_{ g_\kappa } \l_\kappa - \l_\kappa \Ric(g_\kappa)
& = & g_\kappa \\
\l_\kappa |_{ \Sigma} & = & 0 ,
\end{array}
\ri.
\end{equation}
and the metric $ g_\kappa $ is a critical point of the volume functional
$ V(\cdot) $ on the manifold
$ \ \mathcal{M}^{- n (n -1) \kappa^2}_{ \gamma_\kappa} $,
where $  \gamma_\kappa =  g_\kappa |_{T(\Sigma)} $.
For any $ h \in T_{ g_\kappa } \mathcal{M}^{- n (n -1) \kappa^2}_{ \gamma_\kappa} $, define the second variational functional of the volume
\begin{equation}
\begin{split}
\mathcal{F}_\kappa ( h )  = &
\int_B \lf\{  \frac{1}{4}  ( \tr_{g_\kappa} h )^2
 +   \l_\kappa  \lf[  | \div_{g_\kappa} ( h )
 - \frac{1}{2} d ( \tr_{g_\kappa} h ) |_{g_\kappa}^2
+ \frac{1}{4}  | \nabla_\kappa h |_{g_\kappa}^2  \ri]  \ri\}  d V_\kappa \\
& + \int_B  \l_\kappa \lf[ ( \div_{g_\kappa} h)_{l;p} h^{lp}
-   \la h ,  \nabla^2_{g_\kappa} ( \tr_{g_\kappa} h )
 \ra_{g_\kappa} \ri] d V_\kappa  \\
 &  - \int_B  \frac{1}{2} \l_\kappa
  | \div_{g_\kappa} h |_{g_\kappa}^2  d V_\kappa +
\int_B \frac{\kappa^2 }{2} \l_\kappa
\lf[ \lf( \tr_{g_\kappa} h \ri)^2 - | h |^2_{ g_\kappa } \ri] d V_\kappa
\end{split}
\end{equation}
according to (\ref{svar}) in Proposition \ref{svar-formula}, where
$ \nabla_\kappa (\cdot) $ denotes the covariant derivative with
respect to $ g_\kappa$ and $ d V_\kappa $ denotes the volume form
of $ g_\kappa $.

Now, for $ \kappa = 0 $,
choose an $ h_0 \in T_{g_0}  \mathcal{M}^0_{ \gamma_0} $
such that  $ \mathcal{F}_0 ( h_0 ) < 0 $.
The existence of such an $ h_0 $ was proved in (i).
For any $  \kappa \in ( 0,1] $,
by Proposition \ref{deformation-p1},  there is a $ t_\kappa > 0 $
and $ \epsilon_\kappa > 0 $ such that for all $ | t | < t_\kappa $,
$ \hg_\kappa (t) = g_\kappa + t h_0 $ is a smooth metric
on $ B $ and the following Dirichlet boundary value
problem has a unique solution $ u_\kappa (t) $ such that
$ 1 - \epsilon_\kappa \leq u_\kappa (t) \leq 1 + \epsilon_\kappa $:
\begin{equation}
  \begin{cases}
    \alpha\Delta_{\hg_\kappa(t)}u- {R}_\kappa(t)u&=-Ku_\kappa^a,
    \text{\ in $B$}\\
    u_\kappa &=1,  \text{\ on $\Sigma$.}
\end{cases}
    \end{equation}
where $\alpha=4(n-1)/(n-2)$, $a=(n+2)/(n-2)$, and
${R}_\kappa(t)$
is the scalar curvature of $\hg_\kappa(t)$. Moreover,
$v_\kappa = u^\prime_\kappa  (0) $ exists and is a smooth function on $\ol{B}$ which
is the unique solution of:
\begin{equation}\label{vkappaPDE}
  \begin{cases}
    (n-1)\Delta_{g_\kappa} v_\kappa - n ( n -1 ) \kappa^2 v_\kappa & =
    \frac{ n - 1 }{ 4 } {R}^\prime_\kappa ( 0 ),
    \text{\ in $B$}\\
    v_\kappa &=0,  \text{\ on $\Sigma$,}
\end{cases}
    \end{equation}
where
\begin{equation*}
\begin{split}
 {R}^\prime_\kappa  ( 0 ) =
& D \mathcal{R}_{g_\kappa} ( h_0) \\
= &
- \Delta_{g_\kappa} ( \tr_{g_\kappa} h_0 )  +
\div_{g_\kappa} ( \div_{g_\kappa} h_0 )
 - \la h_0 , \Ric(g_\kappa ) \ra_{g_\kappa} .
 \end{split}
\end{equation*}
Since $ g_\kappa \rightarrow g_0 $ in $ C^\infty (\bar{B} )$
as $ \kappa \rightarrow 0 $,
we have $ D \mathcal{R}_{ g_\kappa } ( h_0 ) \rightarrow
 0 $ in  $ C^\infty (\bar{B} )$ and, by (\ref{vkappaPDE}),
 $ v_\kappa \rightarrow 0 $ in $ C^\infty(\bar{B})$.
Now define $ g_\kappa (t)  = u_\kappa^\frac{ 4 }{ n-2 } (t)
\hg_\kappa (t ) $,
$ \{ g_\kappa (t) \}_{ | t | < t_\kappa } $ is a smooth path
in $  \mathcal{M}^{- n (n -1) \kappa^2}_{ \gamma_\kappa} $
with $ g_\kappa ( 0 ) = g_\kappa $.
Let $ h_\kappa = g^\prime_\kappa ( 0 ) $, then
\begin{equation}
h_\kappa  = \frac{ 4 }{ n -2 } v_\kappa g_0 + h_0 .
\end{equation}
Let
$ \kappa \rightarrow 0 $, we have
$ h_\kappa \rightarrow h_0 $ in $ C^\infty(\bar{B}) $ and
$ \mathcal{F}_\kappa ( h_\kappa ) \rightarrow \mathcal{F}_0 (h_0) $.
As $ h_0 $ is chosen so that $ \mathcal{F}_0 ( h_0 ) < 0 $, we see
that there is a small $ \kappa_0 > 0 $, depending only on
$ B $, $ g_0 $ and $ h_0 $, such that
$ \mathcal{F}_\kappa ( h_\kappa )  < 0 $ for $ \kappa < \kappa_0 $.
Recall
\begin{equation*}
\frac{ d^2 }{ d t^2 } V ( g_\kappa ( t ) ) |_{ t = 0 }  =
\mathcal{F}_\kappa (h_\kappa) ,
\end{equation*}
we conclude that
the case $ \Omega \subset \mathbb{H}^n $ is proved by
scaling the metric $ g_\kappa $ to a metric with constant
sectional curvature $ - 1  $.

The case $ \Omega \subset \mathbb{S}^n_+ $ can be
proved in a similar way by replacing $ g_\kappa $
with
\begin{equation}
g_\kappa =  \lf( 1 + \frac{ \kappa^2 }{ 4 } | x |^2 \ri)^{ - 2 } g_0
\end{equation}
and by replacing $ \l_\kappa $ with
\begin{equation}
\l_k = \frac{ 1 }{ ( n-1 ) \kappa^2 }
\lf( - 1 + \frac{ 4 + \kappa^2 }{ 4 - \kappa^2 } \cos{ \kappa r } \ri),
\end{equation}
where $ r $ is the geodesic distance from the origin in
$( B, g_\kappa )$.
This completes the proof of the theorem.
\end{proof}

As a corollary of Theorem \ref{critical-spaceforms-s2-t1},
Theorem \ref{BY-volume-s2-t3} and Theorem \ref{svar-negative},
we have the following nonexistence result on the global volume
minimizer in $ \mathcal{M}^0_\gamma $.

\begin{thm} \label{global-minimum}
Let $\Omega$ be a domain in $\mathbb{R}^n$ bounded by
a smooth, compact, strictly convex hypersurface $ \Sigma $.
Let $ g_0$ be the standard  Euclidean  metric on $ \mathbb{R}^n$
and let $\gamma = g_0 |_{ T(\Sigma) }$.
If $ \Sigma $ is a round sphere in $ \mathbb{R}^n$,
the dimension $ n $ is assumed to satisfy $ 3 \leq n \leq 5$.
Define
$ \beta = \inf \{ V ( g ) \ | \ g \in  \mathcal{M}^0_\gamma \} ,$
then there does not exist a smooth metric $ g $ in
$  \mathcal{M}^0_\gamma  $ such that $ V( g)  = \beta $.
\end{thm}

\begin{proof}
Suppose there exists a smooth metric $ g \in \mathcal{M}^0_\gamma$
with $ V( g)  = \beta $, then $ g $ is a critical point of $V(\cdot)$
in $ \mathcal{M}^0_\gamma $. By Theorem \ref{BY-volume-s2-t3},
$ V(g) \geq V_0 ,$ where $V_0$ is the Euclidean volume of $ \Omega$. Therefore, $ \beta \geq V_0$.

Suppose $ \S $ is not a round sphere in $ \mathbb{R}^n$.
By Theorem \ref{critical-spaceforms-s2-t1},
$ g_0 $ is  not a critical point of $V(\cdot)$ in
$ \mathcal{M}^0_\gamma $. In particular, there is
a path of metrics $ \{g(t) \}$ in $ \mathcal{M}^0_\gamma $
with $ g(0) = g_0 $ such that $ \frac{d}{dt} V(g(t)) |_{t=0} <  0 $.
Hence,  $ V(g(t))  < V_0 $ for small positive $t$,
contradicting $ \beta \geq V_0 $.

If $ \S $ is a round sphere in $ \mathbb{R}^n $, by
Theorem \ref{svar-negative},  there is a path of metrics $ \{g(t) \}$
in ${\mathcal{M}^0_\gamma }$ with $ g(0) = g_0 $
such that $ \frac{d}{dt} V(g(t)) |_{t=0} = 0 $ and
$  \frac{d^2}{d t^2} V(g(t)) |_{t=0} < 0 $. Hence,
$ V(g(t))  < V_0 $ for small $t$, again
contradicting $ \beta \geq V_0 $.
The theorem is proved.
\end{proof}

Before we end this section, we give a discussion of
``large" geodesic balls in $ \mathbb{S}^n $ which strictly
contains  a hemisphere.

\begin{prop} \label{Large-ball}
Let $ \Omega $ be a  geodesic ball  in  $ \mathbb{S}^n$
with geodesic radius $ R $ satisfying 
$ \frac{\pi}{2}  < R < \pi$.
Let $ g $ be the standard metric on $ \mathbb{S}^n$.
Let $ \Sigma $ be the boundary of $ \Omega $ and
 $ \gamma = g |_{ T (\Sigma) }$.
Suppose $ \{ g(t) \}$ is a smooth path of metrics in
$ \mathcal{M}_\gamma $ with  $ g(0) = g $ 
and $ g(t) \in \mathcal{M}^{n(n-1)}_\gamma$.
Then
$$ \frac{d}{dt} V (g(t)) |_{t=0} = 0 . $$
If $ h = g^\prime(0) $
is nonzero and  satisfies
$$ \div_g h = 0 \ \ \mathrm{and} \ \ \tr_g h = 0 , $$
then
$$ \frac{d^2}{dt^2} V ( g(t) ) |_{t=0}  < 0 . $$

\end{prop}

\begin{proof} We embed $ \mathbb{S}^n $ in
$\mathbb{R}^{n+1} $ as the unit sphere
$$ \mathbb{S}^n = \{ ( x_0, x_1, \ldots, x_n ) \ | \
 x_0^2 + x_1^2 + \ldots + x^2_n = 1 \} . $$
Suppose $ \Omega $ is given by the set
$$ \Omega = \mathbb{S}^n \cap \{ x_n > \frac{1}{a} \} $$
for some constant $ a < -1 $.
Consider the function
$$ \l = \frac{1}{ n - 1 } \lf( a x_n - 1 \ri) . $$
Then $ g $ and $ \l $ satisfy
\begin{equation}
   \left\{%
\begin{array}{ll}
    -(\Delta_g\lambda)g+\nabla^2_g\lambda-\lambda \Ric(g)
    &=g \hbox{\ in $\Omega$ ;} \\
    \lambda &=0 \hbox{\ on $\Sigma$.}
\end{array}%
\right.
\end{equation}
From the proof of Theorem \ref{critical-nec-suff-s1t1}, we know
$\frac{ d}{dt } V ( g(t) ) |_{ t = 0 } = 0  $.  By Theorem
\ref{svar-formula} and the assumption $ \div_g h = 0 $ and
$ \tr_g h = 0$, we  have
\begin{equation}
\begin{split}
\frac{d^2}{dt^2} V ( g(t) ) |_{t = 0}
 =  \frac{1}{4}  \int_\Omega  \l    | \nabla h |_g^2
  +  \frac{ 1 }{2  } \int_\Omega  \l  | h |^2_g  .
\end{split}
\end{equation}
Note that, unlike the case $ \Omega \subset \mathbb{S}^n_+ $
in which $ \l $ is positive, we have  $ \l < 0 $ in this case.
Therefore, $ \frac{d^2}{dt^2} V ( g(t) ) |_{t = 0} < 0 $.
\end{proof}

Comparing Proposition \ref{Large-ball} with the case $ \Omega
\subset \mathbb{S}^n_+$ in Theorem \ref{svar-positive},  we find
there exists a dichotomy between  the variational
properties of $ V (\cdot) $ on  ``small"  geodesic balls
in $\mathbb{S}^n $ strictly contained  in a hemisphere and
``large" geodesic balls in  $\mathbb{S}^n $ strictly containing a
hemisphere.

\section{Appendix: TT tensors with prescribed compact support} \label{TT}
In this appendix, we give a construction of trace free and divergence free
$(0,2)$ symmetric tensors with prescribed compact support on 
a rotationally symmetric manifold. As mentioned earlier,
such tensors are used to construct special metric variations in Section 
\ref{sec-secondvar}. 

Let $ m \geq 2 $ be an integer. Let $ S^{ m } $ be the $ m $-dimensional
sphere with the standard differential structure.
Let $ \gst $ be the standard metric on $ S^{ m } $ with constant
sectional curvature $ + 1 $.
For a given $ R > 0 $, consider the product manifold
\be M = ( 0, R) \times S^m \ee
 with a rotationally symmetric metric
\be \label{appendix-eq-1}
g = \frac{ 1 }{ N^2 } d r^2 + r^2 \gst ,
\ee
where $ r $ denotes the usual coordinate on $ (0, R) $ and
$ N = N ( r ) $ is a given smooth positive function on $ ( 0 , R ) $.

In what follows, we let $ \omega $ denote points in $ S^m $ and
let  $ Y = Y (\omega)$ be a fixed eigenfunction of 
$ \mathbb{S}^m = (S^m, \gst)$ whose eigenvalue 
$ \kappa $ is {\em not} the first eigenvalue of $ \mathbb{S}^m$, i.e. 
\be
\label{appendix-eq-2} \Delta_{\gst} Y + \kappa Y = 0 , \ \kappa >
m .
 \ee

\begin{thm}
For any   $ 0 < r_1 < r_2 < R $, and any given smooth function
$a(r)$ with support in $(r_1,r_2)$,  there exists a smooth
divergence free and trace free $ (0, 2) $ symmetric tensor $ h $
on $( M , g )$ such that the support of $ h $ is contained in $
(r_1, r_2 ) \times S^m $ and $h(\p_r, \p_r) =  a(r) Y$.
\end{thm}

\begin{proof}
We first assume that $ a(r), b(r), c(r), d(r) $ are some smooth functions
of the variable $ r $ alone, which  are to
be determined later. For each $ r $, let $ \S_r $ be the leaf
$\{ r \} \times S^m $ in $ M $. We define a $(0,2) $ symmetric tensor
$ h $ on $ M $  as follows,
\be \label{appendix-eq-3}
h(\p_r, \p_r) =  a(r) Y,
\ee
\be \label{appendix-eq-4}
h(\p_r, \cdot)|_{T (\Sigma_r) } = b(r) d Y ,
\ee
\be \label{appendix-eq-5}
h(\cdot, \cdot)|_{ T(\Sigma_r) } = r^2 [ c(r) \hst Y +
d(r) Y \gst ] ,
\ee
where $ h(\p_r, \cdot)|_{T ( \Sigma_r) } $, $ h(\cdot, \cdot)|_{ T(\Sigma_r) } $
are the restriction of $ h(\p_r, \cdot) $, $ h(\cdot, \cdot)$
to the tangent space of $ \Sigma_r $,
$ d Y $ is the differential of $ Y $ on $ S^m $ and
$ \hst Y $ is the Hessian of $ Y $ on $\mathbb{S}^m $.

Let $ \{ \omega_A \ | \ A = 1, \ldots, m \} $ be a local coordinate
chart on $ S^m $ and $ \p_A = \frac{ \p }{\p \omega_A }$ be the
associated tangent vector. Then
\begin{equation} \label{appendix-eq-6}
\begin{split}
\tr_g h = & g^{rr} h_{rr} + g^{AB} h_{AB} \\
 = & N^2 a(r) Y + c(r) \Delta_{\gst} Y+ m d(r) Y  \\
 = & [ N^2 a(r)  -  \kappa c(r)   + m d(r)   ] Y .
\end{split}
\end{equation}

To compute $ \div_g h $, we let $ \p_n = N \p_r $ be the unit
normal vector to $ \S_r  $. Let $ \{ \Gamma^i_{jk} \} $ be the
Christoffel symbol of $ g $ with respect to the frame field $\{
\p_A, \p_n \}$, where $i,j,k \in \{ 1, \ldots, m, n\}$. Let `$;$'
denote the covariant differentiation with respect to $ g $ and
`$,$' denote the partial differentiation. Direct computations
shows \be \label{appendix-eq-7}
\begin{split}
( \div_g h )_n = & h_{nn;n} + g^{AB} h_{n A; B} \\
= & h_{nn;n} + g^{AB} \lf( h_{nA, B} - h_{iA} \Gamma^i_{nB}
- h_{ni} \Gamma^i_{AB} \ri) \\
=  &  h_{nn;n} + g^{AB} \lf( h_{nA, B} - h_{CA} \Gamma^C_{nB}
- h_{nC} \Gamma^C_{AB} - h_{nn} \Gamma^n_{AB} \ri) \\
= &  h_{nn;n} + g^{AB} \lf( h_{nA, B}
- h_{nC} \Gamma^C_{AB} - h_{CA} \Pi^C_{B} + h_{nn} \Pi_{AB} \ri) \\
= &  h_{nn;n} + \div_{\S_r} [ h(\p_n, \cdot) |_{ T( \S_r )} ]
- \la h |_{T(\S_r)},  \Pi \ra + h_{nn} H ,
\end{split}
\ee
where $ \div_{\S_r} ( \cdot) $ denotes the divergence operator
on $ (\S_r, g|_{T( \S_r )})$,
$ \Pi_{AB} = \la \nabla_{\partial_A} \p_n , \partial_B \ra  $
is the second fundamental form of $ \S_r $
and $ H $ is the mean curvature of $ \S_r $.
Similarly,
\be \label{appendix-eq-8}
\begin{split}
( \div_g h)_A = & h_{nA;n} + g^{BC} h_{AB;C} \\
= & h_{nA;n} + g^{BC} \lf(  h_{AB,C} - h_{iB} \Gamma^i_{AC}
- h_{Ai} \Gamma^i_{BC}   \ri) \\
= & h_{n A;n} + ( \div_{\S_r} [ h|_{T( \S_r) } ] )_A
 + g^{BC} h_{nB} \Pi_{AC} + H h_{n A} .
\end{split}
\ee
Meanwhile, by \eqref{appendix-eq-1}, \eqref{appendix-eq-2}, \eqref{appendix-eq-3},
 \eqref{appendix-eq-4} and \eqref{appendix-eq-5}, we have
\begin{itemize}
\item $ H = N m r^{-1} $, $ \Pi_{AB } = N  r^{-1} g_{AB} $
\item $h_{nA}  = N h(\p_r, \p_A ) = N b(r) Y_{,A} $
\item $ h_{nn} = N^2 h(\p_r, \p_r) = N^2 a(r) Y $
\item $ h_{nn;n} = h_{nn, n} = N \p_r [ N^2 a(r) Y ]= N \p_r [ N^2 a(r) ] Y $
\item $  h_{n A;n} =  h_{n A,n} - h_{n i }\Gamma^i_{A n} - h_{iA} \Gamma^i_{nn}
= N \p_r [ N b(r) Y_{,A} ] - N^2 r^{-1} b(r) Y_{,A} $
\item $ \la h |_{ T( \S_r}), \Pi \ra
= N  r^{-1}  [-  \kappa c(r)  + m d(r) ] Y $
\item $ g^{BC} h_{nB} \Pi_{AC} = N^2 r^{-1} b(r) Y_{,A} $
\end{itemize}
and
\begin{equation}
\div_{\S_r} [ h(\p_n, \cdot) |_{ T( \S_r )} ] =
r^{-2} \div_{\gst} [ N b(r) d Y ] = - r^{-2} N b(r) \kappa Y,
\end{equation}
\begin{equation} \label{appendix-eq-9}
\begin{split}
 \div_{\S_r} [ h|_{T( \S_r)}] = & \div_{g_{\mathbb{S}^m}}
 \lf[ c(r) \nabla^2_{\gst} Y + d (r) Y \gst \ri] \\
 = & c(r)  [ d ( \Delta_{\gst}  Y ) + ( m -1) d Y ] + d (r) d Y  \\
 = & \{ c(r)  [ - \kappa  + ( m -1) ] + d (r) \} d Y ,
 \end{split}
 \end{equation}
where in \eqref{appendix-eq-9} we have also used the fact
$ Ric(\gst) = ( m -1) \gst $.
Therefore, it follows from \eqref{appendix-eq-6}, \eqref{appendix-eq-7}
and \eqref{appendix-eq-8} that $ h $ is trace free and divergence free
if and only if the following system holds
\begin{equation}
\left\{
\begin{array}{cll}
 N^2 a  -  \kappa c   + m d  & = & 0 \\
\p_r ( N^2 a ) - r^{-2}  b \kappa -  r^{-1}  [-  \kappa c  + m d ]
+  m r^{-1} N^2 a  & = & 0 \\
N \p_r ( N b )  +  c  [ - \kappa  + ( m -1) ] + d + m r^{-1} N^2 b & = & 0 ,
\end{array}
\right.
\end{equation}
or equivalently
\begin{equation} \label{appendix-eq-TT}
\left\{
\begin{array}{cll}
 -  \kappa c   + m d  & = & - N^2 a  \\
  r^{-2}  b \kappa  & = &  \p_r ( N^2 a ) +  ( m + 1)  r^{-1} N^2 a  \\
c  [ - \kappa  + ( m -1) ] + d  & = & - N \p_r ( N b )  -  m r^{-1} N^2 b .
\end{array}
\right.
\end{equation}
Now, let $ a = a(r) $ be given as in the Theorem. By the second
equation in \eqref{appendix-eq-TT}, $ b = b(r) $ is then
determined accordingly, furthermore $ b(r) $ is smooth and
has compact support in $ (r_1, r_2) $. With $ a(r) $ and $ b(r) $
given, the first and the third equations in \eqref{appendix-eq-TT}
 become a linear system for $ c= c(r) $
and $ d = d(r) $. As
\begin{equation}
\left|
\begin{array}{cc}
- \kappa  & m \\
- \kappa + (m-1) & 1
\end{array}
\right| = ( m -1) ( \kappa - m )
\end{equation}
and $ \kappa $ is chosen so that $ \kappa > m $,
$c(r) $ and $ d(r) $ are uniquely determined by
 the first and the third equations in \eqref{appendix-eq-TT} and
 they both are smooth and have compact support in $(r_1, r_2)$.
With such a choice of $a(r), b(r), c(r), d(r)$, we conclude that
the $(0, 2)$ symmetric tensor  $ h $, defined by
\eqref{appendix-eq-3}--\eqref{appendix-eq-5},
satisfies all the conditions in the Theorem.
\end{proof}

\end{document}